\pdfoutput=1

\documentclass[12pt]{amsart}
\usepackage[margin=1in]{geometry} 

\usepackage{amsmath,amssymb,amsthm,bm,colonequals,graphicx,mathrsfs,mathtools,microtype,stmaryrd}
\usepackage[shortlabels]{enumitem}
\usepackage[colorinlistoftodos]{todonotes}

\numberwithin{equation}{section}

\theoremstyle{plain}

\newtheorem{theorem}{Theorem}[section] 

\newtheorem{lemma}[theorem]{Lemma} 

\newtheorem{proposition}[theorem]{Proposition} 

\newtheorem{proposition-definition}[theorem]{Proposition-Definition} 

\theoremstyle{definition}

\newtheorem{definition}[theorem]{Definition} 

\newtheorem{hypothesis}[theorem]{Hypothesis}

\theoremstyle{remark}

\newtheorem{remark}[theorem]{Remark} 

\newcommand{\CC}{\mathbb{C}}
\newcommand{\FF}{\mathbb{F}}

\newcommand{\PP}{\mathbb{P}}
\newcommand{\QQ}{\mathbb{Q}}
\newcommand{\RR}{\mathbb{R}}

\newcommand{\ZZ}{\mathbb{Z}}

\newcommand{\abs}[1]{\lvert #1 \rvert}
\newcommand{\card}[1]{\lvert #1 \rvert}
\newcommand{\norm}[1]{\lVert #1 \rVert}
\newcommand{\floor}[1]{\lfloor #1 \rfloor}
\newcommand{\ceil}[1]{\lceil #1 \rceil}

\newcommand{\eps}{\epsilon}

\newcommand{\map}{\operatorname}

\newcommand{\ol}{\overline}

\newcommand{\defeq}{\colonequals}

\newcommand{\maps}{\colon}

\newcommand{\belongs}{\subseteq}

\newcommand{\set}[1]{\{#1\}}

\usepackage[pagebackref=true]{hyperref}
\usepackage[alphabetic,initials]{amsrefs}

\begin{document}

\title{Diagonal cubic forms and the large sieve}
\date{}
\author{Victor Y. Wang}
\address{Fine Hall, 304 Washington Road, Princeton, NJ 08540, USA}
\address{Courant Institute, 251 Mercer Street, New York, NY 10012, USA}
\address{IST Austria, Am Campus 1, 3400 Klosterneuburg, Austria}
\email{vywang@alum.mit.edu}

\begin{abstract}
Let $N(X)$ be the number of integral zeros $(x_1,\dots,x_6)\in [-X,X]^6$ of $\sum_{1\le i\le 6} x_i^3$.
Works of Hooley and Heath-Brown imply $N(X)\ll_\epsilon X^{3+\epsilon}$, if one assumes automorphy and GRH for certain Hasse--Weil $L$-functions.
Assuming instead a natural large sieve inequality, we recover the same bound on $N(X)$.
This is part of a more general statement,
for diagonal cubic forms in $\geq 4$ variables,
where we allow approximations to Hasse--Weil $L$-functions.
\end{abstract}

\maketitle

\setcounter{tocdepth}{3}

\section{Introduction}

Fix an integer $m\geq 4$.
Fix integers $F_1,\dots,F_m\in \ZZ\setminus \set{0}$
and let
\begin{equation*}
{\textstyle F(\bm{x}) \defeq \sum_{1\le i\le m} F_ix_i^3},
\end{equation*}
where $\bm{x}=(x_1,\dots,x_m)$.
We are interested in the behavior, as $X\to \infty$, of the point count
\begin{equation*}
N_F(X) \defeq \card{\{\bm{x}\in \ZZ^m\cap [-X,X]^m: F(\bm{x})=0\}}.
\end{equation*}
Certain varieties, $V_{\bm{c},k}$, play a key role.
For each $\bm{c}=(c_1,\dots,c_m)\in \ZZ^m$ and field $k$, let
\begin{equation*}
{\textstyle V_{\bm{c},k}
\defeq \{(\xi_1,\dots,\xi_m)\in \PP^{m-1}_k:
\sum_{1\le i\le m} F_i\xi_i^3 = \sum_{1\le i\le m} c_i\xi_i = 0\}},
\end{equation*}
where $\PP^{m-1}_k$ is the projective space with coordinates $\xi_1,\dots,\xi_m$ over $k$.

In the special case $F=\sum_{1\le i\le 6} x_i^3$, with $m=6$, we abbreviate $N_F(X)$ to $N(X)$.
In this case, building on \cite{hooley1986HasseWeil}, the papers \cite{hooley_greaves_harman_huxley_1997} and \cite{heath1998circle} each proved
\begin{equation}
\label{near-opt}
N(X)\ll_\eps X^{3+\eps},
\end{equation}
assuming Hypothesis~HW of \cite{hooley1986HasseWeil}*{\S6; \cite{heath1998circle}*{\S4}}
for the Hasse--Weil $L$-function
of each smooth variety $V_{\bm{c},\QQ}$ with $\bm{c}\ne \bm{0}$.
Unconditionally, by \cite{vaughan2020some}*{Theorem~1.2},
\begin{equation*}
N(X)\ll_\eps X^{7/2}/(\log{X})^{5/2-\eps}
\end{equation*}
for $X\geq 2$,
via methods stemming from work such as \cites{vaughan1986waring, hall1988divisors, boklan1993reduction, brudern2010asymptotic}.

Hypothesis~HW
practically amounts to automorphy,
plus the Grand Riemann Hypothesis (GRH).
Automorphy remains open \cite{wang2022thesis}*{Appendix~A}.
Hooley suggests that a zero-density hypothesis would suffice in place of GRH \cite{hooley1986HasseWeil}*{p.~51}.
Following the usual paths laid out in \cite{iwaniec2004analytic}*{Theorem~10.4},
a general such density hypothesis is provable assuming
automorphy, a large sieve inequality, and progress on the Grand Lindel\"{o}f Hypothesis (GLH).

In the present paper, we show that
a large sieve inequality \emph{by itself}
would imply \eqref{near-opt}.
The precise large sieve inequality we need will be stated in \S\ref{SEC:basic-setup-and-full-main-result}, as Hypothesis~\ref{hypo:special}.

\begin{theorem}
\label{thm:special}
Suppose $m\in \{5,6\}$.
Assume Hypothesis~\ref{hypo:special}.
Then
\begin{equation}
\label{goal:special}
N_F(X) \ll_\eps X^{3(m-2)/4+\eps},
\end{equation}
for all reals $X\ge 1$ and $\eps>0$.
\end{theorem}

For $m=6$, the exponent in \eqref{goal:special} matches \eqref{near-opt}.
In \S\ref{SEC:basic-setup-and-full-main-result}, we state a more general result, Theorem~\ref{THM:conditional-diagonal-cubic-form-bounds}, valid for all $m\ge 4$.
Our methods might also apply elsewhere \cite{wang2022thesis}*{\S9.1}.
For instance, \cite{restricted_cubic_moments} explains how one may hope to use the modularity of elliptic curves over $\QQ$ to \emph{unconditionally}
produce an absolute constant $\delta>0$ such that $$\card{\set{a\in \ZZ: 1\leq a\leq A} \setminus \set{x^2+y^3+z^3: x,y,z\in \ZZ_{\geq 0}}} \ll A^{6/7-\delta}.$$
This would then improve on the existing bound $O_\eps(A^{6/7+\eps})$ due to Br\"{u}dern \cite{brudern1991ternary}.

\subsection*{Conventions}

We let $\ZZ_{\geq c}\defeq \set{n\in \ZZ: n\geq c}$.
We let $\bm{1}_E\defeq 1$ if a statement $E$ holds, and $\bm{1}_E\defeq 0$ otherwise.
For integers $n\geq 1$, we let $\mu(n)$ denote the M\"{o}bius function.

We write $f\ll_S g$, or $g\gg_S f$, to mean $\abs{f} \leq Cg$ for some $C = C(S)>0$.
The implied constant $C$ is always allowed to depend on $m$ and $F$, in addition to $S$.
We let $O_S(g)$ denote a quantity that is $\ll_S g$.
We write $f\asymp_S g$ if $f\ll_S g\ll_S f$.

\section{Framework and results}
\label{SEC:basic-setup-and-full-main-result}

Let $\mathfrak{D} \defeq 3(\prod_{1\le i\le m} F_i)^{2^{m-2}} \in \ZZ$.
For each $\bm{c}\in \ZZ^m$, let
\begin{equation}
\label{EQN:permissible-discriminant-Delta}
\Delta(\bm{c}) \defeq \mathfrak{D}\,
\prod_{(\upsilon_2,\dots,\upsilon_m)\in \{1,-1\}^{m-1}}\,
\biggl((c_1^3/F_1)^{1/2} + \sum_{2\le i\le m} \upsilon_i(c_i^3/F_i)^{1/2}\biggr)
\in \ZZ.
\end{equation}
For each field $k$ in which $\Delta(\bm{c})$ is invertible,
the variety $V_{\bm{c},k}$ is a smooth complete intersection,
by the Jacobian criterion for smoothness.
Let
\begin{equation}
\label{define-S,SC}
\mathcal{S}\defeq \set{\bm{c}\in \ZZ^m: \Delta(\bm{c})\neq 0},
\qquad \mathcal{S}(C)\defeq \mathcal{S}\cap [-C,C]^m.
\end{equation}

For each $\bm{c}\in \mathcal{S}$ and prime $p$,
we define a local Euler factor $L_p(s, \bm{c})$,
following Serre \cite{serre1969facteurs} and Kahn \cite{kahn2020zeta}*{\S5.6}.
First, choose a prime $\ell\neq p$, and let
$$M(\bm{c},\ell)\defeq H^{m-3}(V_{\bm{c},\ol{\QQ}},\QQ_\ell)/H^{m-3}(\PP^{m-1}_{\ol{\QQ}},\QQ_\ell),$$
where $H^i(W,\QQ_\ell)$ denotes the $i$th $\ell$-adic cohomology group of $W$.
Let $M(\bm{c},\ell)^{I_p} \subseteq M(\bm{c},\ell)$ denote the group of inertia invariants of $M(\bm{c},\ell)$.
Let $\alpha_{\bm{c},j}(p)\in \CC$, for $1\le j\le \dim{M(\bm{c},\ell)^{I_p}}$,
be the geometric Frobenius eigenvalues on $M(\bm{c},\ell)^{I_p}$.
Finally, let
\begin{equation}
\label{define-local-L-factor}
\tilde{\alpha}_{\bm{c},j}(p)
\defeq \frac{\alpha_{\bm{c},j}(p)}{p^{(m-3)/2}},
\qquad
L_p(s, \bm{c})
\defeq \prod_{1\leq j\leq \dim{M(\bm{c},\ell)^{I_p}}}
(1-\tilde{\alpha}_{\bm{c},j}(p)p^{-s})^{-1}.
\end{equation}

On multiplying over $p$, we obtain
for each $\bm{c}\in \mathcal{S}$
a global \emph{Hasse--Weil $L$-function}
\begin{equation}
\label{HW}
L(s, \bm{c})\defeq \prod_p L_p(s, \bm{c})
= \sum_{n\ge 1} \lambda_{\bm{c}}(n) n^{-s},
\end{equation}
for some coefficients $\lambda_{\bm{c}}(n)\in \CC$
defined by expanding the product over $p$.
We now state Hypothesis~\ref{hypo:special}.
It asserts a large sieve inequality, \eqref{ineq:special}, in a certain range.

\begin{hypothesis}
\label{hypo:special}
For all reals $C,N,\eps>0$ with $N\le C^3$,
we have
\begin{equation}
\label{ineq:special}
\sum_{\bm{c}\in \mathcal{S}(C)}\,
\Bigl\lvert
\sum_{n\leq N} v_n
\, \lambda_{\bm{c}}(n)
\Bigr\rvert^2
\ll_{\eps} C^\eps \max(C^m, N)\, \sum_{n\leq N} \abs{v_n}^2
\end{equation}
for all vectors $(v_n)_{1\leq n\leq N}\in \CC^{\floor{N}}$.
\end{hypothesis}

We now make some general comments on $L(s,\bm{c})$.
By \cite{kahn2020zeta}*{\S5.6.3, \S5.6.4} and \cite{laskar2017local}*{Corollary~1.2},
the factors $L_p(s, \bm{c})$ are independent of the choice of $\ell$,
and we have
\begin{equation}
\label{GRC}
\abs{\tilde{\alpha}_{\bm{c},j}(p)} \le 1.
\end{equation}
By \eqref{GRC}, the product and series in \eqref{HW} converge absolutely for $\Re(s)>1$.

We have
$\dim{M(\bm{c},\ell)^{I_p}}
\le \dim{M(\bm{c},\ell)}
\ll_m 1$ by \cite{katz2001sums}*{Corollary of Theorem~3}.
Therefore, by \eqref{GRC}, we have $\lambda_{\bm{c}}(n) \ll_\eps n^\eps$ for all $n\ge 1$.
Thus \eqref{ineq:special} is the large sieve inequality that one would naturally expect to hold.
In fact, \eqref{ineq:special} could potentially hold in the range $N\le C^A$ for any constant $A>0$.
However, we will only need it in the range $N\le C^3$.

The coefficients $\lambda_{\bm{c}}(n)$ can be interpreted geometrically,
but it would take us too far afield to detail anything but the simplest case.
For each $\bm{c}\in \ZZ^m$ and prime $p$, let
\begin{equation*}
E_{\bm{c}}(p) \defeq \frac{\card{\set{\bm{x}\in \FF_p^m:
F(\bm{x})=\bm{c}\cdot\bm{x}=0}} - p^{m-2}}{p-1},
\qquad E^\natural_{\bm{c}}(p) \defeq \frac{E_{\bm{c}}(p)}{p^{(m-3)/2}},
\end{equation*}
where $\bm{c}\cdot\bm{x}\defeq \sum_{1\le i\le m} c_ix_i$.
If $p\nmid \Delta(\bm{c})$, then $M(\bm{c},\ell)^{I_p}=M(\bm{c},\ell)$ and
\begin{equation}
\label{LTF}
\lambda_{\bm{c}}(p)
= \sum_{1\leq j\leq \dim{M(\bm{c},\ell)}} \tilde{\alpha}_{\bm{c},j}(p)
= (-1)^{m-3} E^\natural_{\bm{c}}(p),
\end{equation}
by \eqref{define-local-L-factor} and the Grothendieck--Lefschetz trace formula.

We emphasize that our $L$-functions are normalized differently than in \cites{hooley1986HasseWeil,heath1998circle}.
If $H(s,\bm{c})$ is the $L$-function associated to $V_{\bm{c},\QQ}$ in \cite{heath1998circle}*{\S4},
then $$H(s + \tfrac{m-3}{2}, \bm{c}) = L(s,\bm{c}).$$

\subsection*{Proof framework}

We will analyze $N_F(X)$ using the \emph{delta method}, due to \cites{duke1993bounds,heath1996new}.
This method features some complete exponential sums that we now recall.
Let
\begin{equation}
\label{S_c(n)}
S_{\bm{c}}(n)
\defeq \sum_{\substack{1\leq a\leq n:\\ \gcd(a,n)=1}}\,
\sum_{1\leq x_1,\dots,x_m\leq n} e^{2\pi i(aF(\bm{x}) + \bm{c}\cdot\bm{x})/n},
\qquad S^\natural_{\bm{c}}(n)
\defeq \frac{S_{\bm{c}}(n)}{n^{(m+1)/2}},
\end{equation}
for all $\bm{c}\in \ZZ^m$ and integers $n\geq 1$.
It is known that $S_{\bm{c}}(n)$ is \emph{multiplicative} in $n$,
meaning that $S_{\bm{c}}(1)=1$ and $S_{\bm{c}}(n_1n_2) = S_{\bm{c}}(n_1)S_{\bm{c}}(n_2)$ whenever $\gcd(n_1,n_2)=1$ \cite{heath1998circle}*{Lemma~4.1}.
Thus $S^\natural_{\bm{c}}(n)$ is also multiplicative in $n$.
For each $\bm{c}\in \ZZ^m$, let
\begin{equation*}
\Phi(\bm{c},s)
\defeq \sum_{n\geq1} S^\natural_{\bm{c}}(n) n^{-s}
= \prod_p \Phi_p(\bm{c},s),
\end{equation*} 
where $\Phi_p(\bm{c},s) \defeq \sum_{l\ge 0} S^\natural_{\bm{c}}(p^l) p^{-ls}$.
Ultimately, we will see that $S^\natural_{\bm{c}}(n)$ is related to $\lambda_{\bm{c}}(n)$ in a way that allows us to apply a large sieve inequality, like \eqref{ineq:special}, to the delta method.

Before proceeding, we recall two basic definitions from the theory of Dirichlet series.
For $f, g\maps \ZZ_{\geq 1}\to \CC$,
the \emph{Dirichlet convolution} $f\ast g\maps \ZZ_{\geq 1}\to \CC$ is defined by the formula $$(f\ast g)(n) \defeq \sum_{ab=n} f(a)g(b).$$
A Dirichlet series $\sum_{n\ge 1} f(n)n^{-s}$ is said to be \emph{invertible} if $f(1)\ne 0$,
or equivalently, if there exists $g\maps \ZZ_{\geq 1}\to \CC$ with $(f\ast g)(n) = \bm{1}_{n=1}$.

Our work is based on approximations of Dirichlet series.
For each $\bm{c}\in \mathcal{S}$, let $\Psi(\bm{c},s)$ be an invertible Dirichlet series.
The function $\bm{c}\mapsto \Psi(\bm{c},s)$, from $\mathcal{S}$ to the set of Dirichlet series, will be denoted simply by $\Psi$.
For each $\bm{c}\in \mathcal{S}$ and integer $n\ge 1$, let
\[ b_{\bm{c}}(n),\, a_{\bm{c}}(n),\, a'_{\bm{c}}(n) \]
be the $n^{-s}$ coefficients of the Dirichlet series
\[ \Psi(\bm{c},s),\, \Psi(\bm{c},s)^{-1},\, \Phi(\bm{c},s)/\Psi(\bm{c},s), \]
respectively.
In terms of Dirichlet convolution, this means that
\begin{equation}
\label{baa'-convolution-relations}
(a_{\bm{c}}\ast b_{\bm{c}})(n)=\bm{1}_{n=1},
\qquad
a'_{\bm{c}}=S^\natural_{\bm{c}}\ast a_{\bm{c}},
\qquad
S^\natural_{\bm{c}}=a'_{\bm{c}}\ast b_{\bm{c}}.
\end{equation}
For us, the following particular definition of \emph{approximation} will be convenient.

\begin{definition}
\label{DEFN:approximation-of-Phi}
Call $\Psi$ an \emph{approximation of $\Phi$}
if the following three conditions hold:
\begin{enumerate}
\item If $\bm{c}\in \mathcal{S}$,
then $b_{\bm{c}}(n)$ is multiplicative in $n$.

\item For all $\bm{c}\in \mathcal{S}$,
integers $n\ge 1$,
and reals $\eps>0$,
we have $$\max(\abs{b_{\bm{c}}(n)},\abs{a'_{\bm{c}}(n)})
\ll_\eps n^\eps\sum_{d\mid n}\abs{S^\natural_{\bm{c}}(d)}.$$

\item For all $\bm{c}\in \mathcal{S}$ and primes $p\nmid \Delta(\bm{c})$,
we have $a'_{\bm{c}}(p) \ll p^{-1/2}$.
\end{enumerate}
\end{definition}

\begin{theorem}
\label{THM:example-approximations}
Suppose that for each $\bm{c}\in \mathcal{S}$, we have
\begin{equation*}
\Psi(\bm{c},s)
\in \left\{\Phi(\bm{c},s),\quad
\prod_{p\nmid \Delta(\bm{c})}\Phi_p(\bm{c},s),\quad
\prod_{p\nmid \Delta(\bm{c})} L_p(s, \bm{c})^{(-1)^{m-3}},\quad
L(s, \bm{c})^{(-1)^{m-3}}\right\}.
\end{equation*}
Then $\Psi$ is an approximation of $\Phi$.
\end{theorem}

Theorem~\ref{THM:example-approximations} provides natural examples of approximations.
It will not be used until \S\ref{SEC:finish-proof}, so we defer the proof to that section.
For the rest of \S\ref{SEC:basic-setup-and-full-main-result}, fix an approximation $\Psi$ of $\Phi$.

\subsection*{Hypotheses}

Our main general result, Theorem~\ref{THM:conditional-diagonal-cubic-form-bounds}, will assume that \emph{either} of two specific hypotheses holds.
Our first hypothesis is the following:

\begin{hypothesis}
\label{HYPO:second-moment-hypothesis-for-Psi_1}
For all reals $C, N, \eps > 0$ with $N\leq C^3$, we have
\begin{equation}
\label{INEQ:second-moment-hypothesis}
\sum_{\bm{c}\in \mathcal{S}(C)}\,
\Bigl\lvert
\sum_{n\in I} b_{\bm{c}}(n)
\Bigr\rvert^2
\ll_{\eps} C^\eps
\max(C^m, N)\, N
\end{equation}
for all real intervals $I\belongs (0, N]$.
\end{hypothesis}

The following two remarks may help to clarify the nature of this hypothesis.
\begin{enumerate}
\item If $\Psi = L(s,\bm{c})^{-1}$,
then Hypothesis~\ref{HYPO:second-moment-hypothesis-for-Psi_1} would easily follow from GRH.
On the other hand, if $\Psi = L(s,\bm{c})$,
then Hypothesis~\ref{HYPO:second-moment-hypothesis-for-Psi_1} would follow from GLH plus a technical bound on
$\card{\set{\bm{c}\in \mathcal{S}(C): L(s, \bm{c})\;\textnormal{has a pole at}\;s=1}}$.

\item A \emph{density bound}, namely
$\card{\set{\bm{c}\in \mathcal{S}(C): \abs{\sum_{n\in I} b_{\bm{c}}(n)} \ge N^\sigma}}
\ll_\eps C^{m+\eps}/N^{2\sigma-1}$ for $N\le C^3$ and $\sigma\ge 1/2$,
would follow from Hypothesis~\ref{HYPO:second-moment-hypothesis-for-Psi_1}.
But $C^{m+\eps}/N^{2\sigma-1}$ could be quite large even if $N=C^3$ and $\sigma=1$.
This is unlike in some density applications, e.g.~\cite{iwaniec2004analytic}*{Theorem~10.5}, where further input may be needed near $\sigma=1$.
\end{enumerate}

If $\Psi = L(s,\bm{c})^{-1}$, then Hypothesis~\ref{HYPO:second-moment-hypothesis-for-Psi_1} is perhaps unattractive in that $b_{\bm{c}}(n)$ involves the M\"{o}bius function $\mu(n)$.
We might thus wish to pass from $b_{\bm{c}}(n)$ to $a_{\bm{c}}(n)$.
This is possible, to some extent, in the situation of the following definition:

\begin{definition}
\label{DEFN:standard-Psi_1}
Call $\Psi$ \emph{standard} if
for all $\bm{c}\in \mathcal{S}$,
integers $n\ge 1$,
and reals $\eps>0$,
we have $$\max(\abs{b_{\bm{c}}(n)},\abs{a_{\bm{c}}(n)})
\ll_\eps n^\eps.$$
\end{definition}

Let $\vartheta\in \{0,1\}$ if $\Psi$ is standard,
and let $\vartheta \defeq 0$ if $\Psi$ is non-standard.
Let
\begin{equation}
\label{gamctheta}
\gamma_{\bm{c}}(n)
\defeq (1-\vartheta)\cdot b_{\bm{c}}(n)
+ \vartheta\cdot \mu(n)^2a_{\bm{c}}(n).
\end{equation}
We now come to our main hypothesis:
a large sieve inequality for $\gamma_{\bm{c}}$, in a certain range.

\begin{hypothesis}
\label{HYPO:large-sieve-hypothesis-for-Psi_1-or-1/Psi_1}
For all reals $C, N, \eps > 0$ with $N\leq C^3$, we have
\begin{equation}
\label{INEQ:large-sieve-hypothesis}
\sum_{\bm{c}\in \mathcal{S}(C)}\,
\Bigl\lvert
\sum_{n\leq N} v_n
\, \gamma_{\bm{c}}(n)
\Bigr\rvert^2
\ll_{\eps} C^\eps \max(C^m, N)\, \sum_{n\leq N} \abs{v_n}^2
\end{equation}
for all vectors $(v_n)_{1\leq n\leq N}\in \CC^{\floor{N}}$.
\end{hypothesis}

Again, some brief remarks may be helpful.
\begin{enumerate}
\item When $\vartheta = 1$, the factor $\mu(n)^2$ in \eqref{gamctheta} simply restricts us to square-free moduli $n$.

\item Hypothesis~\ref{HYPO:large-sieve-hypothesis-for-Psi_1-or-1/Psi_1} remains open in general \cite{wang2022thesis}*{Remark~4.1.10}.
\end{enumerate}

\subsection*{Results}

Fix a smooth, compactly supported function $w\maps \RR^m\to \RR$.
Assume that
\begin{equation}
\label{COND:w-supported-away-from-origin}
\bm{0}\notin \ol{\set{\bm{x}\in \RR^m: w(\bm{x})\neq 0}}.
\end{equation}
For reals $X\geq 1$, let
\begin{equation}
\label{EQN:define-N_w(X)}
N_{F,w}(X)
\defeq \sum_{\bm{x}\in \ZZ^m}
w(\bm{x}/X)\, \bm{1}_{F(\bm{x})=0}.
\end{equation}
If $m\ge 5$, then let $N'_{F,w}(X) \defeq N_{F,w}(X)$.
If $m=4$,
then let $\Upsilon$ denote the set of $2$-dimensional rational vector spaces $L$ with $F\vert_{L} = 0$,
and let
\begin{equation}
\label{define-N'(X)}
N'_{F,w}(X)
\defeq \sum_{\bm{x}\in \ZZ^m\setminus (\bigcup_{L\in \Upsilon} L)}
w(\bm{x}/X)\, \bm{1}_{F(\bm{x})=0}.
\end{equation}

\begin{theorem}
\label{THM:conditional-diagonal-cubic-form-bounds}
Assume Hypothesis~\ref{HYPO:large-sieve-hypothesis-for-Psi_1-or-1/Psi_1}
or Hypothesis~\ref{HYPO:second-moment-hypothesis-for-Psi_1}.
Then for some constant $\mathfrak{c}(F,w)\in \RR$, we have
\begin{equation}
\label{goal}
N'_{F,w}(X) - \mathfrak{c}(F,w) X^{m-3}
\ll_\eps X^{3(m-2)/4+\eps},
\end{equation}
for all reals $X\geq 1$ and $\eps>0$.
\end{theorem}

Note that $m,F,w$ are fixed.
In other words, the implied constant in \eqref{goal} is allowed to depend on $m,F,w$ in addition to $\eps$.
Also, for numerical reference,
\begin{equation*}
3(m-2)/4
= 1.5\cdot \bm{1}_{m=4} + 2.25\cdot \bm{1}_{m=5} + 3\cdot \bm{1}_{m=6} + \cdots.
\end{equation*}
In particular, if $5\le m\le 6$,
then $m-3 \le 3(m-2)/4$, and \eqref{goal} simply says
$$N_{F,w}(X) \ll_\eps X^{3(m-2)/4+\eps}.$$

The rest of the paper is devoted to the proof of Theorems~\ref{thm:special}, \ref{THM:example-approximations}, and~\ref{THM:conditional-diagonal-cubic-form-bounds}.
In \S\ref{SEC:conversion}, we reduce
Hypothesis~\ref{HYPO:second-moment-hypothesis-for-Psi_1}
to Hypothesis~\ref{HYPO:large-sieve-hypothesis-for-Psi_1-or-1/Psi_1}.
In \S\S\ref{SEC:delta-method-ingredients}--\ref{SEC:singular-hyperplane-section-contributions},
we recall the delta method for $N_{F,w}(X)$,
then analyze parts of it unconditionally
and parts of it using Hypothesis~\ref{HYPO:second-moment-hypothesis-for-Psi_1}.
In \S\ref{SEC:finish-proof}, we tie together the previous sections to complete the proofs.

\section{A conversion between standard coefficients}
\label{SEC:conversion}

In this section, we prove a useful consequence of Hypothesis~\ref{HYPO:large-sieve-hypothesis-for-Psi_1-or-1/Psi_1}.
First, we record some standard lemmas that will be repeatedly used throughout the paper.

\begin{lemma}
\label{count-h-full-integers}
Let $N, h\in \ZZ_{\ge 1}$.
Then there are at most $O_h(N^{1/h})$ integers $n\in [N, 2N)$ such that $v_p(n)\geq h$ holds for all primes $p\mid n$.
\end{lemma}

\begin{proof}
This is classical; see e.g.~\cite{bateman1958theorem}.
\end{proof}

To proceed, we need to introduce some notation.
We write $u\mid v^\infty$ if there exists $k\in \ZZ_{\ge 1}$ with $u\mid v^k$.
For an integer $c\neq 0$,
we let $\map{sq}(c)$ (resp.~$\map{cub}(c)$)
denote the largest square-full (resp.~cube-full)
positive integer divisor of $c$.
We also let $\map{sq}(0)\defeq 0$.

\begin{lemma}
\label{count-Rinfty-divisors}
Let $N, R\in \ZZ_{\ge 1}$.
Then there are at most $O_\eps(N^\eps R^\eps)$ positive integers $n\le N$ with $n\mid R^\infty$.
\end{lemma}

\begin{proof}
We have $\sum_{n\mid R^\infty} \bm{1}_{n\le N}
\le \sum_{n\mid R^\infty} (N/n)^\eps
= N^\eps \prod_{p\mid R} (1 - p^{-\eps})^{-1}
\ll_\eps N^\eps R^\eps$.
\end{proof}

\begin{lemma}
\label{standard-average-bounds}
Let $N\in \ZZ_{\ge 1}$.
Then the following hold:
\begin{enumerate}
\item We have $$\sum_{n\leq N:\, n=\map{sq}(n)} n^{-1/2}
\ll_\eps N^\eps.$$

\item We have $$\sum_{\abs{c}\le N} \map{sq}(c)^{1/2}
\ll_\eps N^{1+\eps}.$$

\item For any $t\in \RR$, we have
$$\sum_{1\le n\le N} \map{cub}(n)^t
\ll_{t,\eps} N^\eps \max(N, N^{1/3+t}).$$
\end{enumerate}
\end{lemma}

\begin{proof}
(1):
By the $h=2$ case of Lemma~\ref{count-h-full-integers} in dyadic intervals $n\in [2^k, 2^{k+1})$, we have
\begin{equation*}
\sum_{n\leq N:\, n=\map{sq}(n)} n^{-1/2}
\ll \sum_{0\le k\le \log_2{N}} (2^k)^{1/2} (2^k)^{-1/2}
\ll_\eps N^\eps.
\end{equation*}

(2):
There are at most $N/d$ positive integers $n\le N$ with $\map{sq}(n)=d$.
Therefore,
\begin{equation*}
\sum_{\abs{c}\le N} \map{sq}(c)^{1/2}
= 2\sum_{1\le n\le N} \map{sq}(n)^{1/2}
\le 2\sum_{d\leq N:\, d=\map{sq}(d)} \frac{N}{d}\cdot d^{1/2}
\ll_\eps N^{1+\eps},
\end{equation*}
where the last inequality follows from (1).

(3):
There are at most $N/n_3$ positive integers $n\le N$ with $\map{cub}(n)=n_3$.
Thus
\begin{equation*}
\begin{split}
\sum_{1\le n\le N} \map{cub}(n)^t
&\le \sum_{n_3\le N:\, n_3=\map{cub}(n_3)} \frac{N}{n_3}\cdot n_3^t \\
&\ll_t \sum_{0\le k\le \log_2{N}} (2^k)^{1/3} (2^k)^{t-1} N
\ll_{t,\eps} N^\eps \max(N, N^{1/3+t}),
\end{split}
\end{equation*}
by the $h=3$ case of Lemma~\ref{count-h-full-integers} in dyadic intervals $n_3\in [2^k,2^{k+1})$.
\end{proof}

\begin{proposition}
\label{PROP:LSH-implies-SMH}
Fix an approximation $\Psi$ of $\Phi$.
Assume Hypothesis~\ref{HYPO:large-sieve-hypothesis-for-Psi_1-or-1/Psi_1}.
Then Hypothesis~\ref{HYPO:second-moment-hypothesis-for-Psi_1} holds.
\end{proposition}

\begin{proof}
First suppose $\vartheta=0$.
Then $\gamma_{\bm{c}} = b_{\bm{c}}$ by \eqref{gamctheta}.
For $C,N,I$ as in Hypothesis~\ref{HYPO:second-moment-hypothesis-for-Psi_1}, the bound \eqref{INEQ:large-sieve-hypothesis} with $v_n\defeq \bm{1}_{n\in I}$ thus trivially implies \eqref{INEQ:second-moment-hypothesis}, as desired.

Now suppose $\vartheta=1$.
Then in particular, $\Psi$ is standard.
For the rest of the proof, let $C,d,N$ denote \emph{positive} variables.
For integers $d$ and intervals $I$, let
$$A_{\bm{c}}(d,I)\defeq
\sum_{n\in I} \bm{1}_{\gcd(d,n)=1}\,
\mu(n)a_{\bm{c}}(n).$$
We have $\gamma_{\bm{c}}(n) = \mu(n)^2a_{\bm{c}}(n)$ by \eqref{gamctheta}.
Taking $v_n\defeq \bm{1}_{n\in I} \bm{1}_{\gcd(d,n)=1}\,\mu(n)$ in \eqref{INEQ:large-sieve-hypothesis},
and observing that $\mu(n)^3 = \mu(n)$,
we find that Hypothesis~\ref{HYPO:large-sieve-hypothesis-for-Psi_1-or-1/Psi_1} implies
\begin{equation}
\label{appl-LS}
\sum_{\bm{c}\in \mathcal{S}(C)}
\abs{A_{\bm{c}}(d,I)}^2
\ll_{\eps} C^\eps \max(C^m, N)\, N
\end{equation}
uniformly over
reals $C$,
integers $d$,
reals $N\leq C^3$,
and real intervals $I\belongs (0, N]$.

To proceed, we rewrite $b_{\bm{c}}(n)$ using multiplicativity.
First, by \eqref{baa'-convolution-relations}, for primes $p$ we have
$$b_{\bm{c}}(p) = -a_{\bm{c}}(p).$$
Furthermore, an integer $n\geq 1$ can be \emph{uniquely} expressed in the form $n_1d$, where $d$ is square-full, $n_1$ is coprime to $d$, and $n_1$ is square-free.
Therefore, for all $n\geq 1$, we have
\begin{equation}
\label{EQN:b-to-a-conversion}
b_{\bm{c}}(n)
= \sum_{n_1d=n}
\bm{1}_{\gcd(d,n_1)=1}\,
\mu(n_1)a_{\bm{c}}(n_1)
\cdot \bm{1}_{d=\map{sq}(d)}\, b_{\bm{c}}(d).
\end{equation}
We note here that $\mu(n_1)$ is supported on square-free integers $n_1$.

Consider a real $C$,
a real $N\leq C^3$,
and a real interval $I\belongs (0, N]$.
Let $B_{\bm{c}}(I)\defeq
\sum_{n\in I} b_{\bm{c}}(n)$.
By \eqref{EQN:b-to-a-conversion}, we have
\begin{equation*}
\begin{split}
B_{\bm{c}}(I)
&= \sum_{n_1d\in I}
\bm{1}_{\gcd(d,n_1)=1}\,
\mu(n_1)a_{\bm{c}}(n_1)
\cdot \bm{1}_{d=\map{sq}(d)}\, b_{\bm{c}}(d) \\
&= \sum_{d\leq N:\, d=\map{sq}(d)} b_{\bm{c}}(d) \cdot A_{\bm{c}}(d,I/d).
\end{split}
\end{equation*}
By the Cauchy--Schwarz inequality over $d$, it follows that
\begin{equation}
\begin{split}
\label{INEQ:Cauchy-over-square-full-d}
\sum_{\bm{c}\in \mathcal{S}(C)} \abs{B_{\bm{c}}(I)}^2
&\leq \sum_{\bm{c}\in \mathcal{S}(C)}\,
\biggl(\,\sum_{d\leq N:\, d=\map{sq}(d)}
\abs{b_{\bm{c}}(d)}^2 d^{-1/2}\biggr)
\biggl(\,\sum_{d\leq N:\, d=\map{sq}(d)} d^{1/2}
\abs{A_{\bm{c}}(d,I/d)}^2\biggr) \\
&\ll_\eps N^\eps \sum_{\bm{c}\in \mathcal{S}(C)}
\sum_{d\leq N:\, d=\map{sq}(d)} d^{1/2}
\abs{A_{\bm{c}}(d,I/d)}^2,
\end{split}
\end{equation}
by Lemma~\ref{standard-average-bounds}(1),
since $b_{\bm{c}}(d) \ll_\eps d^\eps$ by Definition~\ref{DEFN:standard-Psi_1}.
Yet for all integers $d$, we have
$$\sum_{\bm{c}\in \mathcal{S}(C)} \abs{A_{\bm{c}}(d,I/d)}^2 \ll_\eps C^\eps \max(C^m, N/d) \, (N/d)$$
by \eqref{appl-LS},
since $N/d\leq N\leq C^3$ and $I/d\belongs (0, N/d]$.
Plugging this into \eqref{INEQ:Cauchy-over-square-full-d}, we get
\begin{equation*}
\begin{split}
\sum_{\bm{c}\in \mathcal{S}(C)} \abs{B_{\bm{c}}(I)}^2
&\ll_\eps N^\eps
\sum_{d\leq N:\, d=\map{sq}(d)} d^{1/2}
[C^\eps \max(C^m, N/d) \, (N/d)] \\
&\ll_\eps N^{2\eps} C^\eps \max(C^m, N) \, N,
\end{split}
\end{equation*}
where the second inequality follows from Lemma~\ref{standard-average-bounds}(1)
and the trivial bound $\max(C^m, N/d) \le \max(C^m, N)$.
Thus \eqref{INEQ:second-moment-hypothesis} holds,
uniformly over $C,N,I$.
\end{proof}



\section{Delta method ingredients}
\label{SEC:delta-method-ingredients}

Let $X\geq 1$.
Assume \eqref{COND:w-supported-away-from-origin}, i.e.~that $w$ is supported away from $\bm{0}\in \RR^m$.
Such an assumption is implicit in some of the integral estimates in \cites{heath1996new,heath1998circle}.
Set
\begin{equation}
\label{define-Y}
Y\defeq X^{(\deg F)/2} = X^{3/2}.
\end{equation}
Fix $\eps_0\in (0, 10^{-10}]$ and set
\begin{equation}
\label{define-Z}
Z \defeq
Y/X^{1-\eps_0} = X^{1/2 + \eps_0}.
\end{equation}

Let $\varrho_0(x)\defeq \exp(-(1-x^2)^{-1})$ for $\abs{x}<1$,
and $\varrho_0(x)\defeq 0$ for $\abs{x}\ge 1$.
Let
\begin{equation*}
\varrho(x)\defeq \frac{4\varrho_0(4x-3)}{\int_{y\in \RR} \varrho_0(y)\,dy}.
\end{equation*}
For $x>0$ and $y\in\RR$, let
\begin{equation*}
h(x,y) \defeq \sum_{j\geq 1} \frac{1}{xj}
\left(\varrho(xj) - \varrho{\left(\frac{\abs{y}}{xj}\right)}\right).
\end{equation*}
This is precisely the function $h(x,y)$ defined in \cite{heath1996new}*{\S3}.
For $\bm{c}\in \ZZ^m$ and $n>0$, let
\begin{equation*}
I_{\bm{c}}(n) \defeq \int_{\bm{x}\in\RR^m} w(\bm{x}/X)\,
h(n/Y, F(\bm{x})/Y^2)\, e^{-2\pi i(\bm{c}\cdot\bm{x}/n)} \, d\bm{x}.
\end{equation*}
Let $\norm{\bm{c}}\defeq \max_{1\leq i\leq m}(\abs{c_i})$.
We now recall two standard results on the integral $I_{\bm{c}}(n)$.

\begin{proposition}
[\cite{heath1996new}*{par.~1 of \S7}]
\label{PROP:rigorous-modulus-cutoff-Y}
The functions $n\mapsto I_{\bm{c}}(n)$ are supported on a range of the form $n\le M_0(F,w) Y$,
uniformly over $\bm{c}\in\ZZ^m$,
for some constant $M_0(F,w) > 0$.
\end{proposition}

\begin{lemma}
[\cite{heath1998circle}*{(3.9)}]
\label{LEM:c-aspect-I_c(n)-estimates}
If $\norm{\bm{c}}\geq Z$ and $n\geq 1$,
then $I_{\bm{c}}(n)
\ll_{\eps_0,A} \norm{\bm{c}}^{-A}$,
for all $A>0$.
\end{lemma}

Proposition~\ref{PROP:rigorous-modulus-cutoff-Y} and Lemma~\ref{LEM:c-aspect-I_c(n)-estimates},
together with the trivial bound $\abs{S_{\bm{c}}(n)}\leq n^{1+m}$,
imply
\begin{equation}
\label{INEQ:absolute-tail-decay-bound-in-delta-method}
Y^{-2} \sum_{n\geq 1} \sum_{\norm{\bm{c}} > Z}
n^{-m}\abs{S_{\bm{c}}(n)I_{\bm{c}}(n)}
\ll_{\eps_0,A} X^{-A},
\end{equation}
for all $A>0$.
Here $S_{\bm{c}}(n)$ is defined as in \eqref{S_c(n)}.
By \cite{heath1996new}*{Theorem~2, (1.2)},
we have
\begin{equation}
\label{EQN:un-normalized-delta-method}
(1+O_A(Y^{-A}))\, N_{F,w}(X)
= Y^{-2}\sum_{n\geq1}
\sum_{\bm{c}\in\ZZ^m}
n^{-m}S_{\bm{c}}(n)I_{\bm{c}}(n).
\end{equation}
Equivalently, in terms of $S^\natural_{\bm{c}}(n)$, we have
\begin{equation}
\label{EQN:normalized-delta-method}
(1+O_A(X^{-A}))\, N_{F,w}(X)
= X^{-3}\sum_{n\geq1}
\sum_{\bm{c}\in\ZZ^m} n^{(1-m)/2}S^\natural_{\bm{c}}(n)I_{\bm{c}}(n).
\end{equation}
In view of \eqref{INEQ:absolute-tail-decay-bound-in-delta-method},
analyzing $N_{F,w}(X)$ reduces to understanding the quantity
\begin{equation}
\label{EXPR:main-delta-method-quantity}
\Sigma_0 \defeq
X^{-3} \sum_{n\geq1}
\sum_{\bm{c}\in [-Z,Z]^m} n^{(1-m)/2}S^\natural_{\bm{c}}(n)I_{\bm{c}}(n).
\end{equation}
(Here $I_{\bm{c}}(n)
= I_{\bm{c}}(n)\,\bm{1}_{n\le M_0(F,w) Y}$.
But it is more convenient to keep the factor $\bm{1}_{n\le M_0(F,w) Y}$ implicit,
in order to allow for more flexible technique later on.)

We now recall some standard formulas for $S_{\bm{c}}$ at primes $p$ and prime powers $p^l$.

\begin{proposition}
\label{scp}
Say $p\nmid \bm{c}$.
Then $S^\natural_{\bm{c}}(p)
= E^\natural_{\bm{c}}(p) + O(p^{-1/2})$.
\end{proposition}

\begin{proof}
Let
\begin{equation*}
E(p) \defeq \frac{\card{\set{\bm{x}\in \FF_p^m:
F(\bm{x})=0}} - p^{m-1}}{p-1}.
\end{equation*}
By \cite{heath1998circle}*{p.~680}, we have
$S_{\bm{c}}(p) = p^2 E_{\bm{c}}(p) - p E(p)$
and $E(p) \ll p^{(m-2)/2}$.
Thus $$S_{\bm{c}}(p) = p^2 E_{\bm{c}}(p) + O(p^{m/2}).$$
Now divide by $p^{(m+1)/2}$.
\end{proof}

\begin{proposition}
\label{scp2}
Say $p\nmid \Delta(\bm{c})$.
Then $S_{\bm{c}}(p^l) = 0$ for all integers $l\geq2$.
\end{proposition}

\begin{proof}
This follows immediately from \cite{heath1998circle}*{Lemma~4.4}.
\end{proof}

Fix an approximation $\Psi$ of $\Phi$.
Recall the definition of $\mathcal{S}$ from \eqref{define-S,SC}.
For each $\bm{c}\in \mathcal{S}$,
we have $S^\natural_{\bm{c}} = a'_{\bm{c}}\ast b_{\bm{c}}$ by \eqref{baa'-convolution-relations}.
The following result controls the coefficients $a'_{\bm{c}}$ and $b_{\bm{c}}$.

\begin{proposition}
\label{PROP:S_c-as-linear-combination-of-Psi_1(c,s)-coefficients}
Let $\bm{c}\in \mathcal{S}$.
Then $a'_{\bm{c}}(n)$ is multiplicative in $n$.
Moreover, for all primes $p$ and integers $k\geq 1$, we have
\begin{equation*}
\begin{split}
a'_{\bm{c}}(p)\cdot \bm{1}_{p\nmid \Delta(\bm{c})}
&\ll p^{-1/2}, \\
\max(\abs{a'_{\bm{c}}(p^k)},
\abs{b_{\bm{c}}(p^k)})
&\ll_\eps p^{k\eps} +
p^{k\eps} \sum_{d\mid p^k}\abs{S^\natural_{\bm{c}}(d)}
\cdot \bm{1}_{p\mid \Delta(\bm{c})}.
\end{split}
\end{equation*}
\end{proposition}

\begin{proof}
By \eqref{baa'-convolution-relations},
we have $(a_{\bm{c}}\ast b_{\bm{c}})(n)=\bm{1}_{n=1}$
and $a'_{\bm{c}}=S^\natural_{\bm{c}}\ast a_{\bm{c}}$.
Since $b_{\bm{c}},S^\natural_{\bm{c}}$ are multiplicative,
it follows that $a_{\bm{c}},a'_{\bm{c}}$ are too.
It remains to bound $a'_{\bm{c}}(p^k),b_{\bm{c}}(p^k)$.
When $p\mid \Delta(\bm{c})$,
there is nothing to prove,
since condition~(2) in Definition~\ref{DEFN:approximation-of-Phi} already gives what we want.
Now assume $p\nmid \Delta(\bm{c})$.
Then condition~(3) in Definition~\ref{DEFN:approximation-of-Phi} gives $a'_{\bm{c}}(p)\ll p^{-1/2}$.
On the other hand, $E^\natural_{\bm{c}}(p)\ll 1$ by \eqref{LTF} and \eqref{GRC}.
Therefore, condition~(2) in Definition~\ref{DEFN:approximation-of-Phi} gives
\begin{equation*}
b_{\bm{c}}(p^k),a'_{\bm{c}}(p^k)
\ll_\eps p^{k\eps}\sum_{d\mid p^k}\abs{S^\natural_{\bm{c}}(d)}
\ll p^{k\eps},
\end{equation*}
because $S^\natural_{\bm{c}}(p)
= E^\natural_{\bm{c}}(p)+O(p^{-1/2})\ll 1$ by Proposition~\ref{scp}
and $S^\natural_{\bm{c}}(p^l)
\cdot \bm{1}_{l\geq2} = 0$ by Proposition~\ref{scp2}.
This completes the proof.
\end{proof}

Let $\omega(n)$ denote the number of distinct prime factors of $n$.
The following result, which is due to \cites{hooley1986HasseWeil,heath1998circle},
gives a general pointwise bound on $S^\natural_{\bm{c}}(n)$.
\begin{proposition}
\label{PROP:pointwise-bound}
For some constant $A_F>0$, we have
\begin{equation*}
n^{-1/2}\abs{S^\natural_{\bm{c}}(n)}
\le
A_F^{\omega(n)}
\prod_{1\leq i\leq m} \gcd\bigl(\map{cub}(n)^2, \gcd(\map{cub}(n),\map{sq}(c_i))^3\bigr)^{1/12}
\end{equation*}
for all $\bm{c}\in \ZZ^m$ and integers $n\geq 1$.
\end{proposition}

\begin{proof}
By definition, $S^\natural_{\bm{c}}(n) = n^{-(m+1)/2}S_{\bm{c}}(n)$.
Moreover, since $F$ is diagonal, we have
\begin{equation*}
S_{\bm{c}}(p^l)
\ll_F
p^{l(1+m/2)}
\prod_{1\leq i\leq m} \gcd\bigl(\map{cub}(p^l)^2,\gcd(\map{cub}(p^l),\map{sq}(c_i))^3\bigr)^{1/12},
\end{equation*}
by \cite{heath1998circle}*{(5.1) and (5.2)} for $l\geq 2$ and \cite{heath1983cubic}*{Lemma~11} for $l=1$.
The desired result follows immediately from the multiplicativity of $S_{\bm{c}}$.
\end{proof}

We have stated Proposition~\ref{PROP:pointwise-bound} uniformly over $\bm{c}\in \ZZ^m$.
We proceed to analyze the vectors $\bm{c}$ in sets based on which coordinates $c_i$ are nonzero.
For the rest of \S\ref{SEC:delta-method-ingredients},
we fix a set
\begin{equation}
\label{choose-I}
\mathcal{I}\belongs \set{1,2,\dots,m}.
\end{equation}
Let
\begin{equation}
\label{define-R}
\mathcal{R}\defeq \{\bm{c}\in \ZZ^m\cap [-Z,Z]^m:
\bm{1}_{c_i\ne 0} = \bm{1}_{i\in \mathcal{I}}
\textnormal{ for all $i\in \set{1,2,\dots,m}$}\}.
\end{equation}
By definition, if $\bm{c}\in \mathcal{R}$, then $c_i\ne 0$ if and only if $i\in \mathcal{I}$.

Proposition~\ref{PROP:pointwise-bound} implies that
for all $\bm{c}\in \mathcal{R}$ and integers $n\geq 1$,
we have
\begin{equation}
\label{INEQ:convenient-pointwise-bound}
n^{-1/2}S^\natural_{\bm{c}}(n)
\ll_\eps n^\eps
\map{cub}(n)^{(m-\card{\mathcal{I}})/6}
\prod_{i\in\mathcal{I}} \gcd(\map{cub}(n),\map{sq}(c_i))^{1/4}.
\end{equation}
We will repeatedly use \eqref{INEQ:convenient-pointwise-bound} later in the present paper.
We now turn to $I_{\bm{c}}(n)$.

\begin{lemma}
[\cites{heath1996new,heath1998circle}]
\label{LEM:n-aspect-I_c(n)-estimates}
Assume $\card{\mathcal{I}}\geq 1$.
Then uniformly over $\bm{c}\in \mathcal{R}$,
reals $n\ge 1$,
and integers $k\in \set{0, 1}$,
we have
\begin{equation*}
n^k (\partial/\partial n)^k I_{\bm{c}}(n)
\ll_{k,\eps}
X^{m+\eps}
\left(\frac{X\norm{\bm{c}}}{n}\right)^{1-(m+\card{\mathcal{I}})/4}
\prod_{i\in \mathcal{I}}\left(\frac{\norm{\bm{c}}}{\abs{c_i}}\right)^{1/2}.
\end{equation*}
\end{lemma}

\begin{proof}
By \cite{heath1998circle}*{Lemma~3.2}, since $F$ is diagonal, we have
\begin{equation}
\begin{split}
\label{INEQ:original-n-aspect-integral-bound}
n^k (\partial/\partial n)^k I_{\bm{c}}(n)
&\ll_{k,\eps}
\left(\frac{X\norm{\bm{c}}}{n}\right) X^{m+\eps} \prod_{1\leq i\leq m}
\min\left[
\left(\frac{n}{X\abs{c_i}}\right)^{1/2},
\left(\frac{n}{X\norm{\bm{c}}}\right)^{1/4}
\right] \\
&\le
\left(\frac{X\norm{\bm{c}}}{n}\right) X^{m+\eps}
\prod_{i\in\mathcal{I}} \left(\frac{n}{X\abs{c_i}}\right)^{1/2}
\prod_{i\notin\mathcal{I}} \left(\frac{n}{X\norm{\bm{c}}}\right)^{1/4}.
\end{split}
\end{equation}
After writing $(\frac{n}{X\abs{c_i}})^{1/2} = (\frac{n}{X\norm{\bm{c}}})^{1/2} \, (\frac{\norm{\bm{c}}}{\abs{c_i}})^{1/2}$ in the final line of \eqref{INEQ:original-n-aspect-integral-bound},
the desired inequality follows from the fact that
$1 - \card{\mathcal{I}}/2 - (m-\card{\mathcal{I}})/4 = 1-(m+\card{\mathcal{I}})/4$.
\end{proof}


For later convenience, we now make a definition:
for $\bm{c}\in \ZZ^m$ and integers $N\geq 1$, let
\begin{equation}
\label{EQN:define-convenient-Sobolev-norm-on-I_c(n)}
\norm{I_{\bm{c}}}_{1,\infty;N}
\defeq \sup_{n\in \RR:\, N\leq n\leq 4N}
\left(\abs{I_{\bm{c}}(n)}
+ \abs{n(\partial/\partial n) I_{\bm{c}}(n)}\right).
\end{equation}

In the rest of \S\ref{SEC:delta-method-ingredients}, we will concern ourselves only with $\bm{c}\in \mathcal{R}$ such that $\Delta(\bm{c})\ne 0$.
If $\card{\mathcal{I}}=0$, then no such $\bm{c}$ exist,
because $\mathcal{R}=\{\bm{0}\}$ by \eqref{define-R}.
Therefore, we may and do assume $\card{\mathcal{I}}\ge 1$ for the rest of \S\ref{SEC:delta-method-ingredients}.
To proceed further, we break $\mathcal{R}$ into dyadic pieces.
For each $i\in \mathcal{I}$,
let $C_i\in \set{2^t: t\in \ZZ_{\ge 0}}$ with $1\le C_i\le Z$.
Write
\begin{equation}
\label{dyadic-box}
\mathcal{C}
\defeq \{\bm{c}\in \mathcal{R}:
\abs{c_i}\in [C_i, 2C_i)
\textnormal{ for all $i\in \mathcal{I}$}\},
\qquad
C\defeq \max_{i\in \mathcal{I}}(C_i).
\end{equation}

\begin{proposition}
\label{PROP:dyadic-second-moment-of-absolute-sum-of-error-coefficients-a'_c}
Suppose $N_0\in \ZZ_{\ge 1}$ and $N_0\ll X^{O(1)}$.
Then
\begin{equation*}
\sum_{\bm{c}\in \mathcal{C}:\, \Delta(\bm{c})\neq 0}
\biggl(\,
\sum_{n_0\in [N_0, 2N_0)} \abs{a'_{\bm{c}}(n_0)}\biggr)^{\!2}
\ll_\eps X^\eps N_0^{1+(m-\card{\mathcal{I}})/3} \prod_{i\in\mathcal{I}}C_i.
\end{equation*}
\end{proposition}

\begin{proof}
Consider an integer $n_0\in [N_0, 2N_0)$.
If $n_{\bm{c}}\defeq \prod_{p\mid \Delta(\bm{c})}p^{v_p(n_0)}$
and $n_2\defeq \map{sq}(n_0/n_{\bm{c}})$,
then Proposition~\ref{PROP:S_c-as-linear-combination-of-Psi_1(c,s)-coefficients} implies
\begin{equation*}
\begin{split}
a'_{\bm{c}}(n_0)
&= a'_{\bm{c}}(\tfrac{n_0}{n_{\bm{c}}n_2})
\cdot a'_{\bm{c}}(n_2)
\cdot a'_{\bm{c}}(n_{\bm{c}}) \\
&\ll_\eps (\tfrac{n_0}{n_{\bm{c}}n_2})^{-1/2+\eps}
\cdot n_2^\eps
\cdot \abs{a'_{\bm{c}}(n_{\bm{c}})} \\
&\le n_0^{-1/2+\eps} (n_{\bm{c}}n_2)^{1/2} \abs{a'_{\bm{c}}(n_{\bm{c}})}.
\end{split}
\end{equation*}
Since $n_{\bm{c}}\mid \Delta(\bm{c})^\infty$
and $n_2$ is square-full,
we find, upon summing over $n_0$, that
\begin{equation*}
\begin{split}
\sum_{n_0\in [N_0, 2N_0)} \abs{a'_{\bm{c}}(n_0)}
&\ll_\eps
\sum_{\substack{n_{\bm{c}}n_2 \leq 2N_0: \\ n_{\bm{c}}\mid \Delta(\bm{c})^\infty,\; n_2=\map{sq}(n_2)}}
\frac{N_0}{n_{\bm{c}}n_2}\cdot
N_0^{-1/2+\eps} (n_{\bm{c}}n_2)^{1/2} \abs{a'_{\bm{c}}(n_{\bm{c}})} \\
&\ll_\eps N_0^{1/2+2\eps}
\sum_{\substack{n_{\bm{c}} \leq 2N_0: \\ n_{\bm{c}}\mid \Delta(\bm{c})^\infty}}
n_{\bm{c}}^{-1/2} \abs{a'_{\bm{c}}(n_{\bm{c}})} \\
&\ll_\eps N_0^{1/2+2\eps}(N_0 C)^\eps
\max_{\substack{n_{\bm{c}} \leq 2N_0: \\ n_{\bm{c}}\mid \Delta(\bm{c})^\infty}}
n_{\bm{c}}^{-1/2} \abs{a'_{\bm{c}}(n_{\bm{c}})},
\end{split}
\end{equation*}
where we have used Lemma~\ref{standard-average-bounds}(1) to sum over $n_2 \le 2N_0/n_{\bm{c}}$, and then used Lemma~\ref{count-Rinfty-divisors} to bound the sum over $n_{\bm{c}}$ by a maximum.
Furthermore,
\begin{equation*}
\max_{\substack{n_{\bm{c}} \leq 2N_0: \\ n_{\bm{c}}\mid \Delta(\bm{c})^\infty}}
n_{\bm{c}}^{-1/2} \abs{a'_{\bm{c}}(n_{\bm{c}})}
\ll_\eps N_0^{2\eps}
\max_{\substack{d \leq 2N_0: \\ d\mid \Delta(\bm{c})^\infty}}
d^{-1/2} \abs{S^\natural_{\bm{c}}(d)},
\end{equation*}
since $a'_{\bm{c}}(n_{\bm{c}})
\ll_\eps n_{\bm{c}}^\eps\sum_{d\mid n_{\bm{c}}}\abs{S^\natural_{\bm{c}}(d)}$ by condition~(2) in Definition~\ref{DEFN:approximation-of-Phi}.
But
\begin{equation*}
\begin{split}
\sum_{\bm{c}\in \mathcal{C}:\, \Delta(\bm{c})\neq 0}\,
\max_{\substack{d \leq 2N_0: \\ d\mid \Delta(\bm{c})^\infty}}
d^{-1} \abs{S^\natural_{\bm{c}}(d)}^2
&\leq \sum_{\bm{c}\in \mathcal{C}}
\max_{d \leq 2N_0}
d^{-1} \abs{S^\natural_{\bm{c}}(d)}^2 \\
&\ll_\eps N_0^{(m-\card{\mathcal{I}})/3+2\eps}
\sum_{\bm{c}\in \mathcal{C}}
\prod_{i\in \mathcal{I}} \map{sq}(c_i)^{1/2}
\end{split}
\end{equation*}
by \eqref{INEQ:convenient-pointwise-bound},
since $\gcd(\map{cub}(d),\map{sq}(c_i))^{1/4} \leq \map{sq}(c_i)^{1/4}$.
Yet
\begin{equation}
\label{delicate}
\sum_{\bm{c}\in \mathcal{C}}
\prod_{i\in \mathcal{I}} \map{sq}(c_i)^{1/2}
\ll_\eps \prod_{i\in \mathcal{I}} C_i^{1+\eps},
\end{equation}
by Lemma~\ref{standard-average-bounds}(2).
Proposition~\ref{PROP:dyadic-second-moment-of-absolute-sum-of-error-coefficients-a'_c} follows upon combining the previous four displays.
\end{proof}


We are now prepared to prove a crucial bound for \S\ref{SEC:smooth-contributions}.

\begin{lemma}
\label{LEM:unconditional-second-moment}
Suppose $N_0,N\in \ZZ_{\ge 1}$ and $N_0,N\ll X^{O(1)}$.
Let
\begin{equation*}
Q_{\bm{c}} = \norm{I_{\bm{c}}}_{1,\infty;N}
\, \sum_{n_0\in [N_0, 2N_0)} \abs{a'_{\bm{c}}(n_0)}.
\end{equation*}
Then
\begin{equation*}
\biggl(\,\sum_{\bm{c}\in \mathcal{R}:\, \Delta(\bm{c})\neq 0}
Q_{\bm{c}}^2\biggr)^{\!1/2}
\ll_\eps
X^{m+\eps} N_0^{1/2 + (m-\card{\mathcal{I}})/6} (X/N)^{1-(m+\card{\mathcal{I}})/4}
\max[Z^{1+(\card{\mathcal{I}}-m)/4}, 1].
\end{equation*}
\end{lemma}

\begin{proof}
With notation as in Proposition~\ref{PROP:dyadic-second-moment-of-absolute-sum-of-error-coefficients-a'_c}, consider an element $\bm{c}\in \mathcal{C}$.
Then by \eqref{dyadic-box}, we have
$\abs{c_i} \asymp C_i$ for all $i\in\mathcal{I}$,
whence $\norm{\bm{c}} \asymp C$.
Now \eqref{EQN:define-convenient-Sobolev-norm-on-I_c(n)} and Lemma~\ref{LEM:n-aspect-I_c(n)-estimates} imply
\begin{equation*}
\norm{I_{\bm{c}}}_{1,\infty;N}
\ll_\eps X^{m+\eps} (XC/N)^{1-(m+\card{\mathcal{I}})/4} \prod_{i\in\mathcal{I}}(C/C_i)^{1/2},
\end{equation*}
since $\card{\mathcal{I}}\ge 1$.
By Proposition~\ref{PROP:dyadic-second-moment-of-absolute-sum-of-error-coefficients-a'_c}, it follows that
\begin{equation*}
\sum_{\bm{c}\in \mathcal{C}:\, \Delta(\bm{c})\neq 0}
Q_{\bm{c}}^2
\ll_\eps
X^{2m+3\eps} N_0^{1+(m-\card{\mathcal{I}})/3} (XC/N)^{2-(m+\card{\mathcal{I}})/2} \prod_{i\in \mathcal{I}} C.
\end{equation*}
By \eqref{dyadic-box} we have $1\le C\le Z$,
since $1\le C_i\le Z$ for all $i$.
The quantity $C^{2-(m+\card{\mathcal{I}})/2} \prod_{i\in \mathcal{I}} C = C^{2+(\card{\mathcal{I}}-m)/2}$ is maximized either at $C=Z$ or $C=1$, so we conclude that
\begin{equation*}
\sum_{\bm{c}\in \mathcal{C}:\, \Delta(\bm{c})\neq 0}
Q_{\bm{c}}^2
\ll_\eps
X^{2m+3\eps} N_0^{1+(m-\card{\mathcal{I}})/3} (X/N)^{2-(m+\card{\mathcal{I}})/2} \max[Z^{2+(\card{\mathcal{I}}-m)/2}, 1].
\end{equation*}
Summing over all possibilities for $\mathcal{C}$,
we get
\begin{equation*}
\sum_{\bm{c}\in \mathcal{R}:\, \Delta(\bm{c})\neq 0}
Q_{\bm{c}}^2
\ll_\eps
X^{2m+4\eps} N_0^{1+(m-\card{\mathcal{I}})/3} (X/N)^{2-(m+\card{\mathcal{I}})/2} \max[Z^{2+(\card{\mathcal{I}}-m)/2}, 1].
\end{equation*}
Lemma~\ref{LEM:unconditional-second-moment} follows upon taking a square root.
\end{proof}

Having analyzed $I_{\bm{c}}$ and $a'_{\bm{c}}$ above,
we now concentrate on $b_{\bm{c}}$ for the rest of \S\ref{SEC:delta-method-ingredients}.

\begin{proposition}
\label{PROP:dyadic-second-moment-of-absolute-sum-of-Psi_1-coefficients-b_c}
Let the $C_i$, as well as $\mathcal{C}$ and $C$,
be as specified before Proposition~\ref{PROP:dyadic-second-moment-of-absolute-sum-of-error-coefficients-a'_c}.
Suppose $N_1\in \ZZ_{\ge 1}$ and $N_1\ll X^{O(1)}$.
Then
\begin{equation*}
\sum_{\bm{c}\in \mathcal{C}:\, \Delta(\bm{c})\neq 0}
\biggl(\,
\sum_{n_1\in [N_1, 2N_1]} \abs{b_{\bm{c}}(n_1)}
\biggr)^{\!2}
\ll_\eps X^\eps N_1^{\max(2,1+(m-\card{\mathcal{I}})/3)}\prod_{i\in\mathcal{I}}C_i.
\end{equation*}
\end{proposition}

\begin{proof}
We mimic the proof of Proposition~\ref{PROP:dyadic-second-moment-of-absolute-sum-of-error-coefficients-a'_c}.
Consider an integer $n_1\in [N_1, 2N_1]$.
If $n_{\bm{c}}\defeq \prod_{p\mid \Delta(\bm{c})} p^{v_p(n_1)}$,
then by Proposition~\ref{PROP:S_c-as-linear-combination-of-Psi_1(c,s)-coefficients} and the multiplicativity of $b_{\bm{c}}$,
we have
\begin{equation*}
b_{\bm{c}}(n_1)
= b_{\bm{c}}(n_1/n_{\bm{c}})\,
b_{\bm{c}}(n_{\bm{c}})
\ll_\eps (n_1/n_{\bm{c}})^{\eps}\,
\abs{b_{\bm{c}}(n_{\bm{c}})}
\le n_1^{\eps}\,
\abs{b_{\bm{c}}(n_{\bm{c}})}.
\end{equation*}
Upon summing over $n_1$, then,
\begin{equation*}
\begin{split}
\sum_{n_1\in [N_1, 2N_1]} \abs{b_{\bm{c}}(n_1)}
&\ll_\eps
\sum_{n_{\bm{c}}\leq 2N_1:\, n_{\bm{c}}\mid \Delta(\bm{c})^\infty}
\frac{N_1}{n_{\bm{c}}}\cdot
N_1^{\eps}\,
\abs{b_{\bm{c}}(n_{\bm{c}})} \\
&\ll_\eps N_1^{1+2\eps} C^\eps
\max_{n\leq 2N_1}
n^{-1} \abs{b_{\bm{c}}(n)}
\end{split}
\end{equation*}
by Lemma~\ref{count-Rinfty-divisors}.
Condition~(2) in Definition~\ref{DEFN:approximation-of-Phi} implies
\begin{equation*}
\max_{n\leq 2N_1} n^{-1} \abs{b_{\bm{c}}(n)}
\ll_\eps N_1^{2\eps} \max_{n\leq 2N_1} n^{-1} \abs{S^\natural_{\bm{c}}(n)}.
\end{equation*}
But by \eqref{INEQ:convenient-pointwise-bound}, we have
\begin{equation*}
\sum_{\bm{c}\in \mathcal{C}}
\max_{n\leq 2N_1}
n^{-2} \abs{S^\natural_{\bm{c}}(n)}^2
\ll_\eps N_1^{2\eps}\max(1,N_1^{-1+(m-\card{\mathcal{I}})/3})
\sum_{\bm{c}\in \mathcal{C}}
\prod_{i\in \mathcal{I}} \map{sq}(c_i)^{1/2}.
\end{equation*}
The desired result follows upon combining the last three displays
with \eqref{delicate}.
\end{proof}

\begin{lemma}
\label{LEM:second-moment-of-absolute-sum-of-Psi_1-coefficients-b_c}
Suppose $N_1\in \ZZ_{\ge 1}$ and $N_1\ll X^{O(1)}$.
Then
\begin{equation*}
\sum_{\bm{c}\in \mathcal{R}:\, \Delta(\bm{c})\neq 0}
\biggl(\,
\sum_{n_1\in [N_1, 2N_1]} \abs{b_{\bm{c}}(n_1)}
\biggr)^{\!2}
\ll_\eps X^\eps N_1^{\max(2,1+(m-\card{\mathcal{I}})/3)} Z^{\card{\mathcal{I}}}.
\end{equation*}
\end{lemma}

\begin{proof}
This follows from Proposition~\ref{PROP:dyadic-second-moment-of-absolute-sum-of-Psi_1-coefficients-b_c} upon summing over all possibilities for $\mathcal{C}$.
\end{proof}

We need the following lemma in \S\ref{SEC:smooth-contributions}.
Let
\begin{equation}
\label{DEFN:fixed-beta}
\beta\defeq 1 + 10 \cdot M_0(F,w)
\ll 1.
\end{equation}

\begin{lemma}
Assume Hypothesis~\ref{HYPO:second-moment-hypothesis-for-Psi_1}.
Then
\begin{equation}
\label{INEQ:trivially-improved-conditional-second-moment-bound}
\sum_{\bm{c}\in \mathcal{R}:\, \Delta(\bm{c})\neq 0}\,
\Bigl\lvert
\sum_{n_1\in I}b_{\bm{c}}(n_1)
\Bigr\rvert^2
\ll_\eps \min\left(X^\eps Z^m N_1,
X^\eps Z^{\card{\mathcal{I}}} N_1^{\max(2,1+(m-\card{\mathcal{I}})/3)}\right),
\end{equation}
for all positive integers $N_1\leq \beta Y$ and real intervals $I\belongs [N_1, 2N_1]$.
\end{lemma}

\begin{proof}
The bound $X^\eps Z^m N_1$ in \eqref{INEQ:trivially-improved-conditional-second-moment-bound} follows upon applying \eqref{INEQ:second-moment-hypothesis} with $C=(2\beta)^{1/3} Z$ and $N=2N_1$.
Meanwhile, $X^\eps Z^{\card{\mathcal{I}}}N_1^{\max(2,1+(m-\card{\mathcal{I}})/3)}$ comes from Lemma~\ref{LEM:second-moment-of-absolute-sum-of-Psi_1-coefficients-b_c}.
\end{proof}

\section{Contribution from smooth hyperplane sections}
\label{SEC:smooth-contributions}

Recall the key quantity $\Sigma_0$ from \eqref{EXPR:main-delta-method-quantity}, involving a sum over $\bm{c}\in [-Z,Z]^m$.
In this section, we concentrate on vectors $\bm{c}\in \mathcal{S}(Z) = \mathcal{S}\cap [-Z,Z]^m$, in the notation of \eqref{define-S,SC}.
Let
\begin{equation*}
\Sigma_1 \defeq
X^{-3}\sum_{\bm{c}\in \mathcal{S}(Z)}\,
\sum_{n\geq1}n^{(1-m)/2}S^\natural_{\bm{c}}(n)I_{\bm{c}}(n).
\end{equation*}
We will prove the following result:
\begin{theorem}
Assume Hypothesis~\ref{HYPO:second-moment-hypothesis-for-Psi_1}.
Then
\begin{equation}
\label{INEQ:main-conditional-delta-method-bound}
\Sigma_1
\ll_{\eps_0} X^{3(m-2)/4 + O(\eps_0)}.
\end{equation}
\end{theorem}

For each $n\geq 1$, we have $S^\natural_{\bm{c}}(n) = \sum_{n_0n_1 = n} a'_{\bm{c}}(n_0) b_{\bm{c}}(n_1)$,
since $S^\natural_{\bm{c}}=a'_{\bm{c}}\ast b_{\bm{c}}$ by \eqref{baa'-convolution-relations}.
Thus
\begin{equation}
\label{EXPR:decomposed-main-generic-sum}
\Sigma_1 =
X^{-3} \sum_{\bm{c}\in \mathcal{S}(Z)}\,
\sum_{n_0\geq1} a'_{\bm{c}}(n_0)
\sum_{n_1\geq1} (n_0n_1)^{(1-m)/2}I_{\bm{c}}(n_0n_1)b_{\bm{c}}(n_1).
\end{equation}
By Proposition~\ref{PROP:rigorous-modulus-cutoff-Y},
we have $I_{\bm{c}}(n)=0$ when $n > \beta Y/10$, where $\beta$ is as in \eqref{DEFN:fixed-beta}.
Thus
\begin{equation}
\label{dyadic-Sigma_1-decomp}
\Sigma_1 = X^{-3}
\sum_{\bm{c}\in \mathcal{S}(Z)}\,
\sum_{(N_0,N_1)\in \mathcal{A}}\,
\Diamond_{\bm{c},N_0,N_1},
\end{equation}
where
\begin{equation*}
\begin{split}
\mathcal{A}
&\defeq \set{(N_0,N_1)\in \set{2^t: t\in \ZZ_{\ge 0}}^2:
N_0N_1 \leq \beta Y/10}, \\
\Diamond_{\bm{c},N_0,N_1}
&\defeq \sum_{n_0\in [N_0,2N_0)} a'_{\bm{c}}(n_0)
\sum_{n_1\in [N_1,2N_1)} (n_0n_1)^{(1-m)/2}I_{\bm{c}}(n_0n_1)b_{\bm{c}}(n_1).
\end{split}
\end{equation*}

For convenience,
let $N\defeq N_0N_1$,
let $B_{\bm{c}}(J)\defeq\sum_{n_1\in J}b_{\bm{c}}(n_1)$ for intervals $J$,
and let
\begin{equation*}
\heartsuit_{\bm{c},n_0,N_1}
\defeq \sum_{n_1\in [N_1,2N_1)} (n_0n_1)^{(1-m)/2}I_{\bm{c}}(n_0n_1)b_{\bm{c}}(n_1).
\end{equation*}
Recall $\norm{I_{\bm{c}}}_{1,\infty;N}$ from \eqref{EQN:define-convenient-Sobolev-norm-on-I_c(n)}.
We now have enough notation to state a key lemma:
\begin{lemma}
\label{LEM:uniform-partial-summation}
Let $(N_0,N_1)\in \mathcal{A}$.
Then there exists a probability measure $\nu = \nu_{N_0,N_1}$,
supported on the real interval $[N_1,2N_1]$,
such that for all $\bm{c}\in \mathcal{S}$ and $n_0\in \ZZ\cap [N_0,2N_0)$, we have
\begin{equation}
\label{measure-goal}
\heartsuit_{\bm{c},n_0,N_1}
\ll N^{(1-m)/2}
\, \norm{I_{\bm{c}}}_{1,\infty;N}
\, \int_{x\in [N_1,2N_1]} \abs{B_{\bm{c}}([N_1,x))} \, d\nu(x).
\end{equation}
\end{lemma}

\begin{proof}

Let $\bm{c}\in \mathcal{S}$ and $n_0\in \ZZ\cap [N_0,2N_0)$.
For brevity, let $I(n) = n^{(1-m)/2} I_{\bm{c}}(n)$.
Then
\begin{equation*}
\heartsuit_{\bm{c},n_0,N_1}
= \sum_{n_1\in [N_1,2N_1)} I(n_0n_1)\cdot b_{\bm{c}}(n_1).
\end{equation*}
By partial summation over $n_1$, it follows that
\begin{equation*}
\begin{split}
\abs{\heartsuit_{\bm{c},n_0,N_1}}
&\le \norm{I(r)}_{L^\infty([N,4N])}\,
\abs{B_{\bm{c}}([N_1,2N_1))}
+ n_0\,\norm{I'(r)}_{L^\infty([N,4N])}\,
\sum_{k\in [N_1,2N_1)} \abs{B_{\bm{c}}([N_1,k))} \\
&\ll \norm{I(r)}_{L^\infty([N,4N])}\,
\abs{B_{\bm{c}}([N_1,2N_1))}
+ \frac{N}{N_1}\,\norm{I'(r)}_{L^\infty([N,4N])}\,
\sum_{k\in [N_1,2N_1)} \abs{B_{\bm{c}}([N_1,k))},
\end{split}
\end{equation*}
where $\norm{f(r)}_{L^\infty([N,4N])}\defeq \sup_{r\in [N,4N]} \abs{f(r)}$ for continuous functions $f\maps [N,4N]\to \CC$.
Here
\begin{equation*}
\max(\norm{I(r)}_{L^\infty([N,4N])}, N\, \norm{I'(r)}_{L^\infty([N,4N])})
\ll N^{(1-m)/2} \, \norm{I_{\bm{c}}}_{1,\infty;N}
\end{equation*}
by \eqref{EQN:define-convenient-Sobolev-norm-on-I_c(n)}.
Finally, let $$\nu \defeq \frac12\biggl(\delta_{2N_1}
+ \frac{1}{N_1}\sum_{k\in [N_1,2N_1)} \delta_k\biggr),$$
where $\delta_k$ is the Dirac measure supported on the singleton set $\{k\}$.
Then $\nu$ is a probability measure supported on $[N_1,2N_1]$.
Also, the last three displays imply \eqref{measure-goal}.
\end{proof}

Let $(N_0,N_1)\in \mathcal{A}$.
Let $\mathcal{I}$ and $\mathcal{R}$ be as in \eqref{choose-I} and \eqref{define-R}, respectively.
Since we are presently only interested in $\bm{c}\in \mathcal{S}$, we may and do assume $\card{\mathcal{I}}\geq 1$.
For each $\bm{c}\in \mathcal{S}$, we have
\begin{equation*}
\begin{split}
\abs{\Diamond_{\bm{c},N_0,N_1}}
&\leq \sum_{n_0\in [N_0, 2N_0)}\abs{a'_{\bm{c}}(n_0)\, \heartsuit_{\bm{c},n_0,N_1}} \\
&\ll N^{(1-m)/2}
\, \norm{I_{\bm{c}}}_{1,\infty;N}
\, \sum_{n_0\in [N_0, 2N_0)} \abs{a'_{\bm{c}}(n_0)}\,
\int_{x\in [N_1,2N_1]} \abs{B_{\bm{c}}([N_1,x))} \, d\nu(x),
\end{split}
\end{equation*}
where the first and second inequality are justified by
the triangle inequality and Lemma~\ref{LEM:uniform-partial-summation},
respectively.
Abbreviating $B_{\bm{c}}([N_1,x))$ to $B_{\bm{c}}(x)$ for convenience,
we deduce that
\begin{equation}
\label{INEQ:main-generic-deleted-dyadic-piece}
\sum_{\bm{c}\in \mathcal{R}:\, \Delta(\bm{c})\neq 0}
\abs{\Diamond_{\bm{c},N_0,N_1}}
\ll_\eps X^{m+\eps} Q_1 \, \biggl(\,\sum_{\bm{c}\in \mathcal{R}:\, \Delta(\bm{c})\neq 0}
\left(\int_{x\in [N_1,2N_1]} \abs{B_{\bm{c}}(x)} \, d\nu\right)^{\!2}\,\biggr)^{\!1/2}
\end{equation}
by the Cauchy--Schwarz inequality and Lemma~\ref{LEM:unconditional-second-moment}, where
\begin{equation}
\label{EXPR:main-dyadic-bad-factor-bound-in-delta-method}
Q_1 \defeq N^{(1-m)/2}
N_0^{1/2 + (m-\card{\mathcal{I}})/6}
(X/N)^{1-(m+\card{\mathcal{I}})/4}
\max[Z^{1+(\card{\mathcal{I}}-m)/4}, 1].
\end{equation}
Now, for the rest of \S\ref{SEC:smooth-contributions},
we assume Hypothesis~\ref{HYPO:second-moment-hypothesis-for-Psi_1}.
We have
$$\left(\int_{x\in [N_1,2N_1]} \abs{B_{\bm{c}}(x)} \, d\nu \right)^{\!2}
\ll \int_{x\in [N_1,2N_1]} \abs{B_{\bm{c}}(x)}^2 \, d\nu$$
by the Cauchy--Schwarz inequality,
so
\begin{equation*}
\begin{split}
\biggl(\,\sum_{\bm{c}\in \mathcal{R}:\, \Delta(\bm{c})\neq 0}
\left(\int_{x\in [N_1,2N_1]} \abs{B_{\bm{c}}(x)} \, d\nu\right)^{\!2}\,\biggr)^{\!1/2}
&\ll \biggl(\,\int_{x\in [N_1,2N_1]} \sum_{\bm{c}\in \mathcal{R}:\, \Delta(\bm{c})\neq 0}
\abs{B_{\bm{c}}(x)}^2 \, d\nu\biggr)^{\!1/2} \\
&\ll_\eps X^\eps Q_2
\end{split}
\end{equation*}
by \eqref{INEQ:trivially-improved-conditional-second-moment-bound}, where
\begin{equation}
\label{EXPR:main-dyadic-good-factor-bound-in-delta-method}
Q_2 \defeq \min\left(Z^mN_1,Z^{\card{\mathcal{I}}}N_1^{\max(2,1+(m-\card{\mathcal{I}})/3)}\right)^{1/2}.
\end{equation}

\begin{lemma}
\label{LEM:final-conditional-numerics}
We have $Q_1Q_2 \ll_{\eps_0} X^{3/2-m/4+O(\eps_0)}$.
\end{lemma}

\begin{proof}
We split the proof into four cases.

\emph{Case~1: $\card{\mathcal{I}}=m$.}
Then $Q_2 = (Z^mN_1)^{1/2}$, since $\card{\mathcal{I}}=m$ and $N_1\ge 1$.
Therefore, $Q_1Q_2=Q_3$, where
\begin{equation}
\label{EXPR:generic-dyadic-bound-for-Case-1}
Q_3 \defeq
Z^{m/2}N_1^{1/2}\cdot N^{(1-m)/2}N_0^{1/2}
(X/N)^{1-m/2}\max[Z, 1].
\end{equation}
But $Q_3 = Z^{m/2}X^{1-m/2}\max[Z,1]$, since $N_1N_0 = N$.
By \eqref{define-Z} we have $Z = X^{1/2+\eps_0} \ge 1$,
so $$Q_3 = X^{1-m/2}Z^{1+m/2} = X^{3/2-m/4+(1+m/2)\eps_0}.$$
Thus $Q_1Q_2=Q_3 \ll_{\eps_0} X^{3/2-m/4+O(\eps_0)}$.

\emph{Case~2: $\card{\mathcal{I}}=m-1$ and $N_1\geq Z$.}
Then $Q_2 = (Z^mN_1)^{1/2}$, by \eqref{EXPR:main-dyadic-good-factor-bound-in-delta-method}.
Therefore, $Q_1Q_2=Q_4$, where
\begin{equation*}
Q_4 \defeq
Z^{m/2}N^{1-m/2}N_0^{(m-\card{\mathcal{I}})/6}
(X/N)^{1-(m+\card{\mathcal{I}})/4}
\max[Z^{1+(\card{\mathcal{I}}-m)/4}, 1],
\end{equation*}
since $N_1N_0 = N$.
Since $(m-\card{\mathcal{I}})/6\geq 0$
and $N_0 = N/N_1\le N/Z$, we have
\begin{equation}
\label{rhs1}
Q_4 \le
Z^{m/2}N^{1-m/2}(N/Z)^{(m-\card{\mathcal{I}})/6}
(X/N)^{1-(m+\card{\mathcal{I}})/4}
\max[Z^{1+(\card{\mathcal{I}}-m)/4}, 1].
\end{equation}
The right-hand side of \eqref{rhs1} is \emph{decreasing} as a function of $N$, because
\begin{equation}
\label{case-2-total-N-exponent}
1-m/2 + (m-\card{\mathcal{I}})/6 - 1 + (m+\card{\mathcal{I}})/4 = (\card{\mathcal{I}}-m)/12 < 0.
\end{equation}
Since $N\ge N_1\ge Z$, it follows that
\begin{equation*}
\begin{split}
Q_4 &\le
Z^{m/2}Z^{1-m/2}(Z/Z)^{(m-\card{\mathcal{I}})/6}
(X/Z)^{1-(m+\card{\mathcal{I}})/4}
\max[Z^{1+(\card{\mathcal{I}}-m)/4}, 1] \\
&\ll_{\eps_0} X^{O(\eps_0)}
X^{1-(m+\card{\mathcal{I}})/8} \max[X^{1/2+(\card{\mathcal{I}}-m)/8}, 1],
\end{split}
\end{equation*}
since $Z = X^{1/2+\eps_0}$.
But $\card{\mathcal{I}}=m-1$, so
\begin{equation}
\label{final-|I|=m-1-numerics}
X^{1-(m+\card{\mathcal{I}})/8} \max[X^{1/2+(\card{\mathcal{I}}-m)/8}, 1]
= \max[X^{3/2-m/4}, X^{9/8-m/4}]
= X^{3/2-m/4}.
\end{equation}
Thus $Q_1Q_2=Q_4 \ll_{\eps_0} X^{3/2-m/4+O(\eps_0)}$.

\emph{Case~3: $1\leq \card{\mathcal{I}}\leq m-2$.}
By \eqref{EXPR:main-dyadic-good-factor-bound-in-delta-method},
we have $Q_2 \le (Z^{\card{\mathcal{I}}}N_1^{\max(2,1+(m-\card{\mathcal{I}})/3)})^{1/2}$.
Since $N_1N_0 = N$,
it follows that $Q_1Q_2\le Q_5$, where
\begin{equation*}
Q_5 \defeq
Z^{\card{\mathcal{I}}/2}N_1^{\max(1/2,(m-\card{\mathcal{I}})/6)}N^{1-m/2}N_0^{(m-\card{\mathcal{I}})/6}
(X/N)^{1-(m+\card{\mathcal{I}})/4}
\max[Z^{1+(\card{\mathcal{I}}-m)/4}, 1].
\end{equation*}
Since $N_0\ge 1$ and $N_1N_0=N$,
we have $N_1^{\max(1/2,(m-\card{\mathcal{I}})/6)}N_0^{(m-\card{\mathcal{I}})/6}
\le N^{\max(1/2,(m-\card{\mathcal{I}})/6)}$.
Thus
\begin{equation}
\label{rhs2}
Q_5 \le
Z^{\card{\mathcal{I}}/2}N^{\max(1/2,(m-\card{\mathcal{I}})/6)}N^{1-m/2}
(X/N)^{1-(m+\card{\mathcal{I}})/4}
\max[Z^{1+(\card{\mathcal{I}}-m)/4}, 1].
\end{equation}
The right-hand side of \eqref{rhs2} is \emph{weakly decreasing} in $N$, because
\begin{equation*}
\begin{split}
\max(1/2,(m-\card{\mathcal{I}})/6) + 1-m/2 - 1 + (m+\card{\mathcal{I}})/4
&= \max(1/2,(m-\card{\mathcal{I}})/6) + (\card{\mathcal{I}}-m)/4 \\
&\le 0,
\end{split}
\end{equation*}
in view of the inequality $\card{\mathcal{I}}-m\leq -2$.
Since $N\ge 1$ and $\card{\mathcal{I}}\le m$, it follows that
\begin{equation*}
\begin{split}
Q_5 &\le Z^{\card{\mathcal{I}}/2}X^{1-(m+\card{\mathcal{I}})/4}\max[Z^{1+(\card{\mathcal{I}}-m)/4}, 1] \\
&\le Z^{\card{\mathcal{I}}/2}X^{1-(m+\card{\mathcal{I}})/4}Z \\
&\ll_{\eps_0} X^{O(\eps_0)} X^{3/2-m/4},
\end{split}
\end{equation*}
since $Z = X^{1/2+\eps_0}$.
Thus $Q_1Q_2\le Q_5 \ll_{\eps_0} X^{3/2-m/4+O(\eps_0)}$.

\emph{Case~4: $\card{\mathcal{I}}=m-1$ and $N_1\leq Z$.}
Arguing as in Case~3, we have $Q_1Q_2 \le Q_5$.
But if we hold $N_1$ constant, and plug $N_0 = N/N_1$ into $Q_5$,
then $Q_5$ is \emph{decreasing} in $N$,
by \eqref{case-2-total-N-exponent}.
Since $N\ge N_1$, it follows that $Q_5\le Q_6$, where
\begin{equation*}
Q_6 \defeq Z^{\card{\mathcal{I}}/2}N_1^{\max(1/2,(m-\card{\mathcal{I}})/6)}N_1^{1-m/2}
(X/N_1)^{1-(m+\card{\mathcal{I}})/4}
\max[Z^{1+(\card{\mathcal{I}}-m)/4}, 1].
\end{equation*}
But $Q_6$ is \emph{increasing} in $N_1$,
because
$$\max(1/2,(m-\card{\mathcal{I}})/6) + 1-m/2 - 1 + (m+\card{\mathcal{I}})/4
= 1/4 > 0,$$
in view of the equality $\card{\mathcal{I}}=m-1$.
Since $N_1\le Z$ and $\card{\mathcal{I}}=m-1$, it follows that
\begin{equation*}
\begin{split}
Q_6 &\le Z^{\card{\mathcal{I}}/2}Z^{1/2}Z^{1-m/2}
(X/Z)^{1-(m+\card{\mathcal{I}})/4}\max[Z^{1+(\card{\mathcal{I}}-m)/4}, 1] \\
&= Z\, (X/Z)^{1-(m+\card{\mathcal{I}})/4} \max[Z^{1+(\card{\mathcal{I}}-m)/4}, 1] \\
&\ll_{\eps_0} X^{O(\eps_0)}
X^{1-(m+\card{\mathcal{I}})/8}\max[X^{1/2+(\card{\mathcal{I}}-m)/8}, 1],
\end{split}
\end{equation*}
since $Z = X^{1/2+\eps_0}$.
But $\card{\mathcal{I}}=m-1$, so it follows from \eqref{final-|I|=m-1-numerics} that $Q_6 \ll_{\eps_0} X^{O(\eps_0)} X^{3/2-m/4}$.
Thus $Q_1Q_2\le Q_5\le Q_6 \ll_{\eps_0} X^{3/2-m/4+O(\eps_0)}$.
\end{proof}

\begin{remark}
Interestingly, the quantity $Q_3$ in \eqref{EXPR:generic-dyadic-bound-for-Case-1} is constant over $(N_0,N_1)\in \mathcal{A}$.
\end{remark}

By Lemma~\ref{LEM:final-conditional-numerics}, the left-hand side of \eqref{INEQ:main-generic-deleted-dyadic-piece} is $\ll_{\eps_0} X^{m+O(\eps_0)} X^{3/2-m/4}$.
Upon summing over $(N_0,N_1)\in \mathcal{A}$ and the set of $2^m - 1$ possible sets $\mathcal{R}$,
it follows from \eqref{dyadic-Sigma_1-decomp} that
$$\Sigma_1 \ll_{\eps_0} X^{-3} X^{m+O(\eps_0)} X^{3/2-m/4}
= X^{3(m-2)/4 + O(\eps_0)}.$$
This yields the desired inequality, \eqref{INEQ:main-conditional-delta-method-bound}.

\section{Contribution from the central terms}
\label{SEC:singular-series-contribution}

Here we address the $\bm{c}=\bm{0}$ contribution to \eqref{EXPR:main-delta-method-quantity},
using the theory of $I_{\bm{0}}(n)$ developed in \cite{heath1996new}.
We roughly follow \cite{heath1996new}*{\S12, par.~2}.
Let
\begin{equation}
\label{Sigma2}
\Sigma_2 \defeq
X^{-3}\sum_{n\geq1}
n^{(1-m)/2}S^\natural_{\bm{0}}(n)I_{\bm{0}}(n).
\end{equation}
We begin with a slight extension of \cite{vaughan1997hardy}*{Lemma~4.9}.

\begin{lemma}
\label{LEM:sum-S_0(n)-trivially}
If $N\geq 1$, then $\sum_{n\in [N, 2N)} n^{-m}\abs{S_{\bm{0}}(n)}\ll_\eps N^{(4-m)/3+\eps}$.
\end{lemma}

\begin{proof}
We have $S^\natural_{\bm{0}}(n)
\ll_\eps n^{1/2+\eps}\map{cub}(n)^{m/6}$ by Proposition~\ref{PROP:pointwise-bound}.
Thus $$n^{-m}S_{\bm{0}}(n)\ll_\eps n^{1-m/2+\eps} \map{cub}(n)^{m/6}.$$
Taking $t=m/6$ in Lemma~\ref{standard-average-bounds}(3), we get
\begin{equation*}
\sum_{n\in [N, 2N)} n^{-m}\abs{S_{\bm{0}}(n)}
\ll_\eps N^{1-m/2+\eps} \max(N, N^{1/3+m/6})
= N^{(4-m)/3+\eps},
\end{equation*}
where we note that $\max(N, N^{1/3+m/6}) = N^{1/3+m/6}$
because $N\ge 1$ and $m\ge 4$.
\end{proof}

Lemma~\ref{LEM:sum-S_0(n)-trivially} implies, in particular, the familiar fact that
the singular series
\begin{equation}
\label{EQ:define-singular-series}
\mathfrak{S}
\defeq \sum_{n\geq1} n^{-m}S_{\bm{0}}(n)
\end{equation}
converges absolutely for $m\geq 5$.
It is also known that the real density
\begin{equation}
\label{EQ:define-real-density}
\sigma_{\infty,w}
\defeq \lim_{\eps\to 0}{(2\eps)^{-1}
\int_{\abs{F(\bm{x})}\leq\eps} w(\bm{x}) \, d\bm{x}}
\end{equation}
exists; see e.g.~\cite{heath1996new}*{Theorem~3}.
Yet for all $n\ll Y$, \cite{heath1996new}*{Lemma~13} implies
\begin{equation}
\label{EQ:I0-approximation}
X^{-m}I_{\bm{0}}(n) = \sigma_{\infty,w} + O_A((n/Y)^A),
\end{equation}
for all $A>0$.
If $m\ge 5$,
then via \eqref{EQ:I0-approximation} with $A=(m-4)/3$,
we get
\begin{equation*}
\begin{split}
&\sum_{n\le M_0(F,w) Y} n^{-m}S_{\bm{0}}(n)X^{-m}I_{\bm{0}}(n) \\
&= \sigma_{\infty,w} \sum_{n\le M_0(F,w) Y} n^{-m}S_{\bm{0}}(n)
+ \sum_{n\le M_0(F,w) Y} O((n/Y)^{(m-4)/3} n^{-m}\abs{S_{\bm{0}}(n)}) \\
&= \sigma_{\infty,w}\mathfrak{S}
+ O_\eps(Y^{(4-m)/3 + \eps}),
\end{split}
\end{equation*}
by Lemma~\ref{LEM:sum-S_0(n)-trivially} and \eqref{EQ:define-singular-series}.
Also, by Proposition~\ref{PROP:rigorous-modulus-cutoff-Y}, we have $I_{\bm{0}}(n) = 0$ for all $n>M_0(F,w) Y$.
Since $n^{-m}S_{\bm{0}}(n)=n^{(1-m)/2}S^\natural_{\bm{0}}(n)$ and $Y=X^{3/2}$, it follows that if $m\geq 5$, then
\begin{equation}
\label{Sigma2-estimate}
\Sigma_2 = X^{m-3}\, [\sigma_{\infty,w}\mathfrak{S}
+ O_\eps(X^{(4-m)/2+\eps})]
= \sigma_{\infty,w}\mathfrak{S}X^{m-3}
+ O_\eps(X^{(m-2)/2+\eps}),
\end{equation}
where $\Sigma_2$ is the quantity defined in \eqref{Sigma2}.
On the other hand, for all $m\geq 4$,
\begin{equation}
\label{Sigma2-bound-m=4}
\Sigma_2
\ll X^{m-3}\sum_{n\le M_0(F,w) Y}
n^{-m}\abs{S_{\bm{0}}(n)}
\ll_\eps X^{m-3+\eps}
\end{equation}
by Proposition~\ref{PROP:rigorous-modulus-cutoff-Y} and Lemma~\ref{LEM:sum-S_0(n)-trivially},
since $I_{\bm{0}}(n) \ll X^m$ by \cite{heath1996new}*{Lemma~16}.

\section{Contribution from singular hyperplane sections}
\label{SEC:singular-hyperplane-section-contributions}

In this section, we study the quantity
\begin{equation}
\label{EXPR:main-singular-delta-method-quantity}
\Sigma_3 \defeq X^{-3} \sum_{n\ge 1}\,
\sum_{\bm{c}\in [-Z,Z]^m:\, \Delta(\bm{c}) = 0,\; \bm{c}\ne \bm{0}}
n^{(1-m)/2}S^\natural_{\bm{c}}(n)I_{\bm{c}}(n).
\end{equation}
We will prove the following result, extending work of Heath-Brown.
Recall the definitions of $N_{F,w}(X)$ and $N'_{F,w}(X)$ from \eqref{EQN:define-N_w(X)} and \eqref{define-N'(X)}, respectively.

\begin{theorem}
If $m\ge 5$, then
\begin{equation}
\label{Sigma3-bound}
\Sigma_3 \ll_{\eps_0} X^{3(m-2)/4+O(\eps_0)}.
\end{equation}
If $m=4$, then
\begin{equation}
\label{Sigma3-estimate-m=4}
\Sigma_3 = N_{F,w}(X)-N'_{F,w}(X) + O_{\eps_0}(X^{3(m-2)/4+O(\eps_0)}).
\end{equation}
\end{theorem}

The cases $m=4$ and $m=6$ of this result are due to Heath-Brown.
For instance,
the estimate \eqref{Sigma3-estimate-m=4} for $m=4$ follows directly from \cite{heath1998circle}*{Lemmas~7.2 and~8.1},
in view of the tail estimate \eqref{INEQ:absolute-tail-decay-bound-in-delta-method}.
Therefore, we may and do assume $m\ge 5$, for the rest of \S\ref{SEC:singular-hyperplane-section-contributions}.

We combine ideas from \cite{hooley1986HasseWeil} and \cite{heath1998circle}.
Let $\mathcal{I}$ and $\mathcal{R}$ be as in \eqref{choose-I} and \eqref{define-R}, respectively.
Since we are only interested in $\bm{c}\ne \bm{0}$, we may and do assume $\card{\mathcal{I}}\geq 1$.
Let $\mathcal{C}$ and $C$ be as in \eqref{dyadic-box},
for some $C_i\in \set{2^t: t\in \ZZ_{\ge 0}}$ with $1\le C_i\le Z$.

By Proposition~\ref{PROP:rigorous-modulus-cutoff-Y}, the sum $\Sigma_3$ from \eqref{EXPR:main-singular-delta-method-quantity} is supported on $n\le M_0(F,w) Y$.
Let
\begin{equation}
\label{EXPR:main-dyadic-singular-quantity}
\Sigma_4\defeq X^{-3}\sum_{n\le M_0(F,w) Y}
\sum_{\bm{c}\in \mathcal{C}:\, \Delta(\bm{c}) = 0}
n^{(1-m)/2}S^\natural_{\bm{c}}(n)I_{\bm{c}}(n).
\end{equation}
Now consider an element $\bm{c}\in\mathcal{C}$ with $\Delta(\bm{c})=0$,
assuming such a $\bm{c}$ exists.
Denote the nonempty fibers of the map
$\mathcal{I}\to \QQ^\times/(\QQ^\times)^2,\;
i\mapsto F_ic_i\bmod{(\QQ^\times)^2}$ by
\begin{equation*}
\mathcal{I}(k) \defeq \set{i\in\mathcal{I}: F_ic_i\equiv g_k\bmod{(\QQ^\times)^2}},
\end{equation*}
for $1\le k\le K$, say,
where the $g_k$ are signed, nonzero square-free integers.
Trivially, we have $\sum_{1\le k\le K} \card{\mathcal{I}(k)} = \card{\mathcal{I}}$.
For each $i\in \mathcal{I}(k)$,
we may write
\begin{equation}
\label{EQN:singular-ci-square-constraint}
c_i = g_kF_i^{-1}e_i^2
\end{equation}
with $e_i\in\ZZ$.
Moreover, by \eqref{EQN:permissible-discriminant-Delta}
and the $\QQ$-linear independence of square roots of distinct square-free integers,
we may choose the signs of the integers $e_i$ so that
\begin{equation}
\label{EQN:singular-ci-equation-constraint}
\sum_{i\in \mathcal{I}(k)} F_i(e_i/F_i)^3 = 0.
\end{equation}
Since $c_i\neq 0$ implies $e_i\neq 0$ for all $i\in \mathcal{I}(k)$,
we immediately deduce from \eqref{EQN:singular-ci-equation-constraint} that
\begin{equation}
\label{INEQ:no-singletons-I(k)}
\card{\mathcal{I}(k)}\ge 2.
\end{equation}

We now prove a general lemma that will allow us, in Lemma~\ref{LEM:fixed-n-S_c-second-moment}, to exploit the structure uncovered in the previous paragraph.

\begin{lemma}
\label{LEM:divisor-sieve-bound-for-ge_i^2}
Let $J\in \ZZ_{\ge 2}$, let $d_1,\dots,d_J\in \ZZ_{\ge 1}$,
and let $G,E_1,\dots,E_J\in \RR_{>0}$.
Then
\begin{equation*}
\sum_{\substack{1\leq g\leq G: \\ \mu(g)^2 = 1}}\,
\prod_{1\leq i\leq J}
\sum_{\substack{1\leq e_i\leq E_i: \\ d_i\mid \map{sq}(ge_i^2)}} d_i^{1/2}
\leq \prod_{1\leq i\leq J} (2^{\omega(d_i)} G^{1/2} E_i).
\end{equation*}
\end{lemma}

\begin{proof}
By H\"{o}lder's inequality over $g$, we may assume that $E_1=\dots=E_J=E$ and $d_1=\dots=d_J=d$, say.
Let $S\defeq \set{h\mid d: \mu(h)^2=1}$.
Now consider integers $g,e\geq 1$ with $g$ square-free.
Then $\map{sq}(g e^2) = \gcd(g,e) e^2$.
Therefore, if $d\mid \map{sq}(ge^2)$,
and we let $h\defeq \gcd(g,e,d)$,
then $$h\in S, \qquad (d/h)\mid e^2,$$
whence $e$ is divisible by the integer $\prod_{p\mid (d/h)} p^{\ceil{v_p(d/h)/2}} \ge (d/h)^{1/2}$.
Thus, given $h\in S$, the number of possible $e\in [1,E]$ is at most $E/(d/h)^{1/2}$.
It follows that
\begin{equation}
\label{INEQ:sum-over-e-and-simplify}
\sum_{\substack{1\leq e\leq E: \\ d\mid \map{sq}(ge^2)}} d^{1/2}
\leq \sum_{h\in S} (d^{1/2}\cdot \bm{1}_{h\mid g}\cdot E/(d/h)^{1/2})
= \sum_{h\in S} (\bm{1}_{h\mid g}\cdot h^{1/2}E),
\end{equation}
for every square-free $g\ge 1$.
By \eqref{INEQ:sum-over-e-and-simplify}, and H\"{o}lder's inequality over $h$, we get
\begin{equation*}
\begin{split}
\sum_{\substack{1\leq g\leq G: \\ \mu(g)^2 = 1}}\,
\biggl(\,\sum_{\substack{1\leq e\leq E: \\ d\mid \map{sq}(ge^2)}} d^{1/2}\biggr)^{\!J}
&\leq \sum_{\substack{1\leq g\leq G: \\ \mu(g)^2 = 1}}\,
\biggl(\,\sum_{h\in S} (\bm{1}_{h\mid g}\cdot h^{1/2}E)\biggr)^{\!J} \\
&\leq \sum_{\substack{1\leq g\leq G: \\ \mu(g)^2 = 1}}\, \card{S}^{J-1}
\sum_{h\in S} (\bm{1}_{h\mid g}\cdot h^{1/2}E)^J \\
&\leq \card{S}^{J-1} \sum_{\substack{h\in S: \\ h\leq G}} (G/h) (h^{1/2} E)^J \\
&\le \card{S}^J G^{J/2} E^J,
\end{split}
\end{equation*}
where in the last step we note that $h^{J/2-1} \le G^{J/2-1}$.
This suffices, since $\card{S} = 2^{\omega(d)}$.
\end{proof}

\begin{lemma}
\label{LEM:fixed-n-S_c-second-moment}
Let $n\geq 1$ be an integer.
Then
\begin{equation*}
\sum_{\bm{c}\in \mathcal{C}:\, \Delta(\bm{c})\neq 0}
n^{-1}\abs{S^\natural_{\bm{c}}(n)}^2
\ll_\eps n^\eps \map{cub}(n)^{(m-\card{\mathcal{I}})/3} \prod_{i\in \mathcal{I}} C_i^{1/2+\eps}.
\end{equation*}
\end{lemma}

\begin{proof}
Let $n_3\defeq \map{cub}(n)$.
Fix a set $\mathcal{J}\belongs \mathcal{I}$ with $\card{\mathcal{J}}\geq 2$.
Let $G \in \set{2^t: t\in \ZZ_{\ge 0}}$, and let $E_i\defeq (2F_iC_i/G)^{1/2}$ for each $i\in\mathcal{J}$.
Let $\tau(\cdot)$ be the divisor function.
Then
\begin{equation*}
\begin{split}
&\sum_{\substack{\abs{g}\in [G,2G): \\ \mu(\abs{g})^2 = 1}}\,
\prod_{i\in\mathcal{J}}
\sum_{\substack{\abs{e_i}\leq (2F_iC_i/\abs{g})^{1/2}: \\ gF_i^{-1}e_i^2\in \ZZ\setminus \set{0}}}
\gcd(n_3,\map{sq}(gF_i^{-1}e_i^2))^{1/2} \\
&\le 2^{1+\card{\mathcal{J}}} \sum_{\substack{g\in [G,2G): \\ \mu(g)^2 = 1}}\,
\prod_{i\in\mathcal{J}} \sum_{1\leq e_i\leq E_i} \gcd(n_3,\map{sq}(ge_i^2))^{1/2} \\
&\leq 2^{1+\card{\mathcal{J}}}\sum_{\substack{g\in [G,2G): \\ \mu(g)^2 = 1}}\,
\prod_{i\in\mathcal{J}} \sum_{d_i\mid n_3} \sum_{\substack{1\leq e_i\leq E_i: \\ d_i\mid \map{sq}(ge_i^2)}} d_i^{1/2} \\
&\le 2^{1+\card{\mathcal{J}}}\prod_{i\in\mathcal{J}} (\tau(n_3)^2 G^{1/2} E_i) \\
&= 2^{1+\card{\mathcal{J}}}\tau(n_3)^{2\card{\mathcal{J}}} \prod_{i\in\mathcal{J}} (2F_iC_i)^{1/2},
\end{split}
\end{equation*}
where in the penultimate step we use Lemma~\ref{LEM:divisor-sieve-bound-for-ge_i^2} for each possible choice of divisors $d_i\mid n_3$,
and we note that $2^{\omega(d_i)} \le 2^{\omega(n_3)} \le \tau(n_3)$.
Moreover, if $G > \min_{i\in \mathcal{J}}(2F_iC_i)$, then
\begin{equation*}
\sum_{\substack{\abs{g}\in [G,2G): \\ \mu(\abs{g})^2 = 1}}\,
\prod_{i\in\mathcal{J}}
\sum_{\substack{\abs{e_i}\leq (2F_iC_i/\abs{g})^{1/2}: \\ gF_i^{-1}e_i^2\in \ZZ\setminus \set{0}}}
\gcd(n_3,\map{sq}(gF_i^{-1}e_i^2))^{1/2}
= 0,
\end{equation*}
since the sum over one of the variables $e_i\in \ZZ\setminus\{0\}$ is empty.
On summing the penultimate display over all $G \in \set{2^t: t\in \ZZ_{\ge 0}}$ with $G \le \min_{i\in \mathcal{J}}(2F_iC_i)$, we conclude that
\begin{equation}
\label{blah}
\sum_{\substack{g\in \ZZ\setminus\{0\}: \\ \mu(\abs{g})^2 = 1}}\,
\prod_{i\in\mathcal{J}}
\sum_{\substack{\abs{e_i}\leq (2F_iC_i/\abs{g})^{1/2}: \\ gF_i^{-1}e_i^2\in \ZZ\setminus \set{0}}}
\gcd(n_3,\map{sq}(gF_i^{-1}e_i^2))^{1/2}
\ll_{\card{\mathcal{J}},\eps} n_3^\eps \prod_{i\in\mathcal{J}} (2F_iC_i)^{1/2+\eps}.
\end{equation}

Recall the constraints \eqref{EQN:singular-ci-square-constraint} and \eqref{INEQ:no-singletons-I(k)} on $\set{\bm{c}\in \mathcal{C}: \Delta(\bm{c}) = 0}$.
Applying \eqref{blah} with $\mathcal{J} = \mathcal{I}(k)$,
for each $k\in [1,K]$,
and multiplying the resulting $K$ inequalities, we get
\begin{equation}
\label{blah2}
\prod_{1\le k\le K}
\sum_{\substack{g_k\in \ZZ\setminus\{0\}: \\ \mu(\abs{g_k})^2 = 1}}\,
\prod_{i\in\mathcal{I}(k)}
\sum_{\substack{\abs{e_i}\leq (2F_iC_i/\abs{g_k})^{1/2}: \\ g_kF_i^{-1}e_i^2\in \ZZ\setminus \set{0}}}\,
\gcd(n_3,\map{sq}(g_kF_i^{-1}e_i^2))^{1/2} \\
\ll_\eps n_3^\eps \prod_{i\in \mathcal{I}} C_i^{1/2+\eps},
\end{equation}
since $K\le \card{\mathcal{I}}\le m$, and the variables $m,F_i$ are fixed.
On summing \eqref{blah2} over all possible choices for the sets $\mathcal{I}(k)\subseteq \mathcal{I}$, we deduce that
\begin{equation}
\label{INEQ:main-technical-singular-gcd-bound}
\sum_{\bm{c}\in \mathcal{C}:\, \Delta(\bm{c}) = 0}
\prod_{i\in\mathcal{I}} \gcd(n_3,\map{sq}(c_i))^{1/2}
\ll_\eps n_3^\eps \prod_{i\in \mathcal{I}} C_i^{1/2+\eps}.
\end{equation}
Lemma~\ref{LEM:fixed-n-S_c-second-moment} follows immediately from \eqref{INEQ:convenient-pointwise-bound} and \eqref{INEQ:main-technical-singular-gcd-bound}.
\end{proof}

\begin{remark}
Interestingly, the proof of \eqref{INEQ:main-technical-singular-gcd-bound} uses the constraint \eqref{EQN:singular-ci-equation-constraint} only through \eqref{INEQ:no-singletons-I(k)}.
\end{remark}

Taking $n_3=1$ in \eqref{INEQ:main-technical-singular-gcd-bound} implies
$$\card{\set{\bm{c}\in \mathcal{C}: \Delta(\bm{c}) = 0}}\ll_\eps \prod_{i\in \mathcal{I}} C_i^{1/2+\eps}.$$
Therefore, Lemma~\ref{LEM:fixed-n-S_c-second-moment} implies
\begin{equation}
\label{INEQ:fixed-n-S_c-first-moment}
\sum_{\bm{c}\in \mathcal{C}:\, \Delta(\bm{c}) = 0} n^{-1/2}\abs{S^\natural_{\bm{c}}(n)}
\ll_\eps n^\eps \map{cub}(n)^{(m-\card{\mathcal{I}})/6} \prod_{i\in \mathcal{I}} C_i^{1/2+\eps},
\end{equation}
by the Cauchy--Schwarz inequality over $\bm{c}$.

Let $N\in \set{2^t: t\in \ZZ_{\ge 0}}$ with $1\le N\le M_0(F,w) Y$.
By Lemma~\ref{LEM:n-aspect-I_c(n)-estimates}, \eqref{INEQ:fixed-n-S_c-first-moment},
and the $t=(m-\card{\mathcal{I}})/6$ case of Lemma~\ref{standard-average-bounds}(3),
the sum
\begin{equation*}
\Sigma_5 \defeq X^{-3}\sum_{n\in [N, 2N)}
\sum_{\bm{c}\in \mathcal{C}:\, \Delta(\bm{c}) = 0}
n^{(1-m)/2}\abs{S^\natural_{\bm{c}}(n)I_{\bm{c}}(n)}
\end{equation*}
satisfies the bound $\Sigma_5 \ll_\eps X^{m-3+\eps} Q_7$, where
\begin{equation*}
\begin{split}
Q_7 &\defeq N^{1-m/2} (XC/N)^{1-(m+\card{\mathcal{I}})/4} \max(N,N^{1/3+(m-\card{\mathcal{I}})/6})\, C^{\card{\mathcal{I}}/2} \\
&= X^{1-(\card{\mathcal{I}}+m)/4} \max(N^{1+(\card{\mathcal{I}}-m)/4}, N^{1/3+(\card{\mathcal{I}}-m)/12})\, C^{1+(\card{\mathcal{I}}-m)/4}.
\end{split}
\end{equation*}
Since $N^{1+(\card{\mathcal{I}}-m)/4} = (N^{1/3+(\card{\mathcal{I}}-m)/12})^3$,
we will analyze $Q_7$ according to the sign of $$\mathfrak{e}\defeq 1+(\card{\mathcal{I}}-m)/4.$$

\emph{Case~1: $\mathfrak{e}\le 0$.}
Then, since $N,C\gg 1$, we have $$Q_7 \ll X^{1-(\card{\mathcal{I}}+m)/4} \le X^{(6-m)/4},$$
where the final inequality holds because $X\ge 1$ and $\card{\mathcal{I}}\geq-2$.

\emph{Case~2: $\mathfrak{e}\ge 0$.}
Then, since $N\ll Y$ and $C\ll Z$, we have
\begin{equation*}
Q_7 \ll X^{1-(\card{\mathcal{I}}+m)/4} Y^{1+(\card{\mathcal{I}}-m)/4} Z^{1+(\card{\mathcal{I}}-m)/4}.
\end{equation*}
Plugging in \eqref{define-Y} and \eqref{define-Z}, we get
$$Q_7 \ll_{\eps_0} X^{1-(\card{\mathcal{I}}+m)/4+O(\eps_0)}(X^2)^{1+(\card{\mathcal{I}}-m)/4}
= X^{3+(\card{\mathcal{I}}-3m)/4+O(\eps_0)}.$$
Moreover, if $\card{\mathcal{I}}\leq 2m-6$, then $3+(\card{\mathcal{I}}-3m)/4 \leq (6-m)/4$.

If $m\geq 6$,
then $1\le \card{\mathcal{I}}\le m\le 2m-6$,
so regardless of what $\card{\mathcal{I}}$ is, it follows that
\begin{equation*}
\begin{split}
\Sigma_5 &\ll_{\eps_0} X^{m-3+\eps_0} Q_7 \\
&\ll_{\eps_0} X^{m-3+\eps_0}X^{(6-m)/4+O(\eps_0)} \\
&= X^{3(m-2)/4+O(\eps_0)},
\end{split}
\end{equation*}
whence by summing over all possibilities for $N$ and $\mathcal{C}$ we get
\begin{equation*}
\Sigma_3,\Sigma_4
\ll_{\eps_0} X^{3(m-2)/4+O(\eps_0)},
\end{equation*}
where $\Sigma_3,\Sigma_4$ are as defined in \eqref{EXPR:main-singular-delta-method-quantity} and \eqref{EXPR:main-dyadic-singular-quantity}, respectively.
This completes the proof of \eqref{Sigma3-bound} for $m\ge 6$.
For the rest of \S\ref{SEC:singular-hyperplane-section-contributions}, we relinquish the previous definitions of $\mathcal{C}$ and $C$.

For $m=5$, we first show that a natural extension of \cite{heath1998circle}*{Lemma~7.1} holds.
\begin{lemma}
\label{analogHBLem7.1}
If $5\le m\le 6$ and $C\gg 1$, then
$\card{\{\bm{c}\in \ZZ^m\cap [-C,C]^m: \Delta(\bm{c})=0\}} \ll_\eps C^{m-3+\eps}$.
\end{lemma}

\begin{proof}
For $m=6$, this follows directly from \cite{heath1998circle}*{Lemma~7.1}.
Now suppose $m=5$.
A \emph{partition} of $m$ is an infinite, weakly decreasing sequence of nonnegative integers $\lambda_1,\lambda_2,\dots$, such that $\sum_{k\ge 1} \lambda_k = m$.
For any partition of $m$, let $$e_k \defeq 2\cdot \bm{1}_{2\le \lambda_k\le 4} + (\lambda_k-2)\cdot \bm{1}_{\lambda_k \ge 5}$$ for $k\ge 1$.
Let $\theta$ denote the maximum value of $\frac12\sum_{k\ge 1} e_k$ over all partitions of $m$.
By \cite{heath1998circle}*{p.~687}, we have
$\card{\{\bm{c}\in \ZZ^m\cap [-C,C]^m: \Delta(\bm{c})=0\}} \ll_\eps C^{\theta+\eps}$.

Clearly $\lambda_3\le \floor{m/3} = 1$, so $e_k=0$ for all $k\ge 3$.
If $\lambda_2\le 1$, then $e_k=0$ for all $k\ge 2$, so $\sum_{k\ge 1} e_k = e_1 \le m-2$.
If $\lambda_2\ge 2$, then $\lambda_1\le m-\lambda_2\le 3$, so $e_k\le 2$ for all $k\ge 1$, whence $\sum_{k\ge 1} e_k = e_1+e_2 \le 4$.
In either case, $\sum_{k\ge 1} e_k \le 4$.
Therefore, $\theta\le 2 = m-3$.
\end{proof}

We now recall a bound from \cite{heath1998circle}
that is valid for all $m\ge 4$.

\begin{lemma}
\label{HBp689}
Fix $\varepsilon>0$.
Suppose $1\ll N\ll X^{3/2}$ and $1\ll C\ll X^{1/2+\varepsilon}$.
Let
\begin{equation*}
A = \sum_{N<q\le 2N}\, \sum_{C<\norm{\bm{c}}\le 2C:\, \Delta(\bm{c}) = 0}\,
q^{-m} S_{\bm{c}}(q) I_{\bm{c}}(q).
\end{equation*}
Then there exist reals $X_1,X_2,X_3\gg 1$ and an integer $H\ge 1$
such that $X_1X_2X_3 \asymp N$ and
\begin{equation*}
A \ll_{\varepsilon} X^{m+4\varepsilon}
N^{-m} X_1^{1+m/2} X_2^{2/3+2m/3} X_3^{1+2m/3}
H^{1/2} \left(\frac{N}{XC}\right)^{\!(m-2)/4} \mathcal{N}_1\mathcal{N}_2(H),
\end{equation*}
where in terms of the quantity $\mathfrak{D} = 3(\prod_{1\le i\le m} F_i)^{2^{m-2}}$ from \S\ref{SEC:basic-setup-and-full-main-result}, we let
\begin{equation*}
\begin{split}
\mathcal{N}_1 &\defeq \sum_{(q_1,q_2,q_3):\, X_i<q_i\le 2X_i}
\bm{1}_{\map{cub}(q_1)=1}
\bm{1}_{q_2=\map{cub}(q_2)}
\bm{1}_{q_3\mid \mathfrak{D}^\infty}, \\
\mathcal{N}_2(H) &\defeq \sum_{C<\norm{\bm{c}}\le 2C}
\bm{1}_{H\mid \bm{c}}
\bm{1}_{\Delta(\bm{c})=0}.
\end{split}
\end{equation*}
\end{lemma}

\begin{proof}
This is immediate from \cite{heath1998circle}*{pp.~688--689,
from the definition of $A$ on p.~688
to the definition of $\mathcal{N}_2(H)$ on p.~689}.
What Heath-Brown calls $P$ (resp.~$X$), we call $X$ (resp.~$N$).
Moreover, in terms of Heath-Brown's notation $n$ and $G$,
our $m$ and $\Delta$ satisfy $m=n$ and $\Delta(\bm{c}) = 3G(\bm{c})$.
However, our $C,q,\bm{c},X_1,X_2,X_3,H$ match Heath-Brown's notation.
\end{proof}

Applying Lemma~\ref{count-h-full-integers} to $q_2$ and Lemma~\ref{count-Rinfty-divisors} to $q_3$, it is clear that $$\mathcal{N}_1 \ll_\varepsilon X_1X_2^{1/3}X_3^\varepsilon.$$
Now assume $5\le m\le 6$.
Then $\mathcal{N}_2(H)=0$ unless $H\le 2C$, in which case $$\mathcal{N}_2(H) \ll_\varepsilon (C/H)^{m-3+\varepsilon}$$ by Lemma~\ref{analogHBLem7.1}.
Plugging the last two displays into Lemma~\ref{HBp689}, with $\varepsilon\defeq \eps_0$, we get
\begin{equation*}
A \ll_{\eps_0} X^{m+O(\eps_0)}
N^{-m} X_1^{2+m/2} X_2^{1+2m/3} X_3^{1+2m/3}
H^{1/2} \left(\frac{N}{XC}\right)^{\!(m-2)/4} \left(\frac{C}{H}\right)^{\!m-3}.
\end{equation*}
Since $m-3\geq 1/2$, we have $H^{1/2} (C/H)^{m-3} \le C^{m-3}$.
Moreover, $m\le 6$ implies $2+m/2 \geq 1+2m/3$,
so $X_1^{2+m/2} X_2^{1+2m/3} X_3^{1+2m/3} \ll (X_1X_2X_3)^{2+m/2} \asymp N^{2+m/2}$.
Thus
\begin{equation}
\label{rhs3}
A \ll_{\eps_0} X^{m+O(\eps_0)}
N^{2-m/2}
\left(\frac{N}{XC}\right)^{\!(m-2)/4} C^{m-3}.
\end{equation}
Since $2-m/2+(m-2)/4 = (6-m)/4\ge 0$ (resp.~since $m-3\ge (m-2)/4$),
the right-hand side of \eqref{rhs3} is weakly increasing in $N$ (resp.~in $C$).
Therefore
\begin{equation*}
A \ll_{\eps_0} X^{m+O(\eps_0)}
(X^{3/2})^{2-m/2}
(X^{1/2})^{m-3}
= X^{3m/4+3/2+O(\eps_0)}.
\end{equation*}
Summing over $1\ll N=M_0(F,w)Y/2^{k_1}$ and $1\ll C=Z/2^{k_2}$ with $k_1,k_2\in \ZZ_{\ge 1}$, we get
\begin{equation*}
\Sigma_3 \ll X^{-3} X^{3m/4+3/2+O(\eps_0)}
= X^{3(m-2)/4+O(\eps_0)},
\end{equation*}
where $\Sigma_3$ is the quantity defined in \eqref{EXPR:main-singular-delta-method-quantity}.
This completes the proof of \eqref{Sigma3-bound}.


\section{Proof of main results}
\label{SEC:finish-proof}

In this section, we first prove Theorem~\ref{THM:conditional-diagonal-cubic-form-bounds}, because it builds directly on our work in \S\S\ref{SEC:delta-method-ingredients}--\ref{SEC:singular-hyperplane-section-contributions} on the delta method.
We then prove Theorem~\ref{THM:example-approximations} using \eqref{GRC}, \eqref{LTF}, and Proposition~\ref{scp}.
Finally, we combine Theorems~\ref{THM:example-approximations} and~\ref{THM:conditional-diagonal-cubic-form-bounds} to prove Theorem~\ref{thm:special}.

\begin{proof}
[Proof of Theorem~\ref{THM:conditional-diagonal-cubic-form-bounds}]

By Proposition~\ref{PROP:LSH-implies-SMH},
we see that Hypothesis~\ref{HYPO:large-sieve-hypothesis-for-Psi_1-or-1/Psi_1} implies Hypothesis~\ref{HYPO:second-moment-hypothesis-for-Psi_1}.
Therefore, we may and do assume Hypothesis~\ref{HYPO:second-moment-hypothesis-for-Psi_1}.
Now recall the quantity $\Sigma_0$ from \eqref{EXPR:main-delta-method-quantity}.
By \eqref{EQN:normalized-delta-method} and the tail estimate \eqref{INEQ:absolute-tail-decay-bound-in-delta-method}, we have
\begin{equation*}
N_{F,w}(X) - \Sigma_0 \ll_{A,\eps_0} X^{-A}.
\end{equation*}

\emph{Case~1: $m=4$.}
Then adding \eqref{INEQ:main-conditional-delta-method-bound}, \eqref{Sigma2-bound-m=4}, and \eqref{Sigma3-estimate-m=4} together, we get
\begin{equation*}
\Sigma_0 = \Sigma_1+\Sigma_2+\Sigma_3
= N_{F,w}(X)-N'_{F,w}(X)
+ O_{\eps_0}(X^{3(m-2)/4+O(\eps_0)})
+ O_{\eps_0}(X^{m-3+\eps_0}).
\end{equation*}
It follows that $N'_{F,w}(X) \ll_{\eps_0} X^{3(m-2)/4+O(\eps_0)}$.
Let $\mathfrak{c}(F,w) \defeq 0$.

\emph{Case~2: $m\ge 5$.}
Then adding \eqref{INEQ:main-conditional-delta-method-bound}, \eqref{Sigma2-estimate}, and \eqref{Sigma3-bound} together, we get
\begin{equation*}
\Sigma_0 = \Sigma_1+\Sigma_2+\Sigma_3
= \mathfrak{c}(F,w)X^{m-3} + O_{\eps_0}(X^{3(m-2)/4+O(\eps_0)})
+ O_{\eps_0}(X^{(m-2)/2+\eps_0}),
\end{equation*}
where $\mathfrak{c}(F,w) \defeq \sigma_{\infty,w}\mathfrak{S}$.
It follows that $N_{F,w}(X)-\mathfrak{c}(F,w)X^{m-3} \ll_{\eps_0} X^{3(m-2)/4+O(\eps_0)}$.

In each case, taking $\eps_0 \to 0$ gives the desired result, \eqref{goal}.
\end{proof}

\begin{proof}
[Proof of Theorem~\ref{THM:example-approximations}]
Let $\bm{c}\in \mathcal{S}$.
Since $\Psi(\bm{c},s)$ has an Euler product, condition~(1) in Definition~\ref{DEFN:approximation-of-Phi} clearly holds.
It remains to prove that conditions~(2) and~(3) hold.

\emph{Case~1: $\Psi(\bm{c},s) = \Phi(\bm{c},s)$.}
Then conditions~(2) and~(3) are trivial,
since $$(b_{\bm{c}}(n), a'_{\bm{c}}(n))
= (S^\natural_{\bm{c}}(n), \bm{1}_{n=1}).$$

\emph{Case~2: $\Psi(\bm{c},s) = \prod_{p\nmid \Delta(\bm{c})}\Phi_p(\bm{c},s)$.}
Then conditions~(2) and~(3) are trivial,
since $$(b_{\bm{c}}(n), a'_{\bm{c}}(n))
= (S^\natural_{\bm{c}}(n)\cdot \bm{1}_{\gcd(n,\Delta(\bm{c}))=1},
S^\natural_{\bm{c}}(n)\cdot \bm{1}_{n\mid \Delta(\bm{c})^\infty}).$$

\emph{Case~3: $\Psi(\bm{c},s) \in \{\prod_{p\nmid \Delta(\bm{c})} L_p(s, \bm{c})^{(-1)^{m-3}},
L(s, \bm{c})^{(-1)^{m-3}}\}$.}
Then by \eqref{GRC}, we have
\begin{equation}
\label{L-is-standard}
b_{\bm{c}}(n),a_{\bm{c}}(n)\ll_\eps n^\eps.
\end{equation}
But $a'_{\bm{c}}=S^\natural_{\bm{c}}\ast a_{\bm{c}}$, by \eqref{baa'-convolution-relations}.
Therefore, condition~(2) holds.
Furthermore,
if $p\nmid \Delta(\bm{c})$,
then $a_{\bm{c}}(p)
= -b_{\bm{c}}(p)$ by \eqref{baa'-convolution-relations}
and $b_{\bm{c}}(p) = (-1)^{m-3} \lambda_{\bm{c}}(p) = E^\natural_{\bm{c}}(p)$ by \eqref{LTF},
so $$a'_{\bm{c}}(p)
= S^\natural_{\bm{c}}(p)+a_{\bm{c}}(p)
= S^\natural_{\bm{c}}(p)-E^\natural_{\bm{c}}(p)
\ll p^{-1/2}$$
by Proposition~\ref{scp}.
Therefore, condition~(3) also holds.
\end{proof}

\begin{proof}
[Proof of Theorem~\ref{thm:special}]
Let $\Psi \defeq L(s,\bm{c})^{(-1)^{m-3}}$.
Then $\Psi$ is an approximation of $\Phi$, by Theorem~\ref{THM:example-approximations}.
Moreover, $\Psi$ is standard by \eqref{L-is-standard} and Definition~\ref{DEFN:standard-Psi_1}.
Now let $\vartheta \defeq 1$.
Then $\gamma_{\bm{c}}(n) = \mu(n)^m \lambda_{\bm{c}}(n)$ by \eqref{gamctheta},
since for all primes $p$ we have $a_{\bm{c}}(p) = (-1)^{m-2} \lambda_{\bm{c}}(p)$ by the definition of $a_{\bm{c}}$.
Upon plugging in $\mu(n)^mv_n$ for $v_n$ in Hypothesis~\ref{hypo:special},
we immediately find that Hypothesis~\ref{HYPO:large-sieve-hypothesis-for-Psi_1-or-1/Psi_1} holds.
Let $\varsigma\maps \RR\to \RR$ be a nonnegative, smooth, compactly supported function such that
$\varsigma(t)=1$ for all $t\in [1,4]$,
and $\varsigma(t)=0$ for all $t\notin [\frac12,8]$.
Let $$w(\bm{x}) \defeq \varsigma({\textstyle \sum_{1\le i\le m} x_i^2}).$$
Then Theorem~\ref{THM:conditional-diagonal-cubic-form-bounds}
implies $N_{F,w}(X) \ll_\eps X^{3(m-2)/4+\eps}$ for all $X\ge 1$.
Since $w(\bm{x}/2^k)=1$ for all $\bm{x}\in \ZZ^m$ in the annulus $4^k\le \sum_{1\le i\le m} x_i^2\le 4^{k+1}$, it follows that
\begin{equation*}
\begin{split}
N_F(X) - 1
&= \card{\set{\bm{x}\in \ZZ^m\cap [-X,X]^m: F(\bm{x})=0,\; \bm{x}\ne \bm{0}}} \\
&\le \sum_{0\le k\le \log_4(4mX^2)} N_{F,w}(2^k) \\
&\ll_\eps \sum_{0\le k\le \log_4(4mX^2)} (2^k)^{3(m-2)/4+\eps} \\
&\ll ((4mX^2)^{1/2})^{3(m-2)/4+\eps} \\
&\ll_\eps X^{3(m-2)/4+\eps},
\end{split}
\end{equation*}
for all $X\ge 1$.
This implies Theorem~\ref{thm:special}.
\end{proof}

\section*{Acknowledgements}

I thank Peter Sarnak
for suggesting projects
that ultimately led to the present paper.{\let\thefootnote\relax\footnote{This work was partially supported by NSF grant DMS-1802211, and the European Union's Horizon~2020 research and innovation program under the Marie Sk\l{}odowska-Curie Grant Agreement No.~101034413.}}
I also thank him for many
encouraging discussions, helpful comments, and references.
Thanks also to Tim Browning, Trevor Wooley, and Nina Zubrilina for helpful comments,
and to Levent Alp\"{o}ge and Will Sawin
for some interesting old discussions.
I thank
Yang Liu, Evan O'Dorney, Ashwin Sah, and Mark Sellke
for conversations illuminating the combinatorics of an older, counting version of the present Lemma~\ref{LEM:unconditional-second-moment}.
Finally, special thanks are due to the editors and referees for their patience and help with the exposition.

\bibliographystyle{amsxport}
\bibliography{final.bib}

\end{document}